\documentclass[11pt]{amsart}
\usepackage{amsfonts,amscd} 
\usepackage[matrix,arrow,frame,curve]{xy}
\newcommand{\andname}{et}
\newcommand{\editorname}{\'editeur}

\newcommand{\volumename}{volume}

\newcommand{\ofname}{de}

\newcommand{\Inname}{Dans}

\newcommand{\pagesname}{pages}

\title{Stuctures du cube et fibr{\'e}s d'intersection}
\author{Fran\c{c}ois Ducrot}
\address{D{\'e}partement de Math{\'e}matiques, Universit{\'e} d'Angers\\
2 boulevard Lavoisier, 49045 Angers cedex 01, France}
\email{francois.ducrot@univ-angers.fr}
\date{}
\swapnumbers\newtheorem{theo}{Theor\`eme}[subsection]
\newtheorem{prop}[theo]{Proposition}
\newtheorem{lemme}[theo]{Lemme}
\newtheorem{cor}[theo]{Corollaire}
\newtheorem{ptheo}{Th{\'e}or{\`e}me principal}[section]
\theoremstyle{definition}
\begingroup
\newtheorem{Def}[theo]{D\'efinition}
\newtheorem{rem}[theo]{Remarque}
\newtheorem{ex}[theo]{Exemple}
\newtheorem{nota}[theo]{Notation}
\newtheorem{constr}[theo]{Construction}
\endgroup
\newcommand{\ga }{\alpha }
\newcommand{\gb }{\beta }
\newcommand{\gc }{\gamma }
\newcommand{\gd }{\delta }
\renewcommand{\ge }{\varepsilon }
\newcommand{\gs}{\sigma}

\renewcommand{\L}{\Lambda}
\newcommand{\s}{\sigma}

\renewcommand{\t}{\theta}

\newcommand{\R}{\mathbb{R}}
\newcommand{\Z}{\mathbb{Z}}
\newcommand{\F}{\longrightarrow}

\newcommand{\isom}{\stackrel{\sim }{\F }}

\newcommand{\oo}{\mathcal{O}}
\newcommand{\cc}{\mathcal{C}}
\newcommand{\dd}{\mathcal{D}}

\newcommand{\ox}{\mathcal{O}_{X}}
\newcommand{\os}{\mathcal{O}_{S}}

\newcommand{\tx}{\otimes _{\ox}}

\newcommand{\ts}{\otimes _{\os}}

\newcommand{\inv}{^{-1}}
\newcommand{\lst}{_{\ast }}
\newcommand{\ust}{^{\ast }}
\newcommand{\ve}{^{\vee }}
\newcommand{\pr}{^{\prime}}
\renewcommand{\sec}{^{\prime \prime}}
\newcommand{\ddp}{\mathcal{D}^{\prime}}
\newcommand{\gdp}{\delta ^{\prime}}
\newcommand{\drpi}{\det \text{R}\pi \lst }
\newcommand{\fisec}{\phi^{\prime \prime\ast}}
\newcommand{\fipr}{\phi^{\prime\ast}}
\begin{document}
\setcounter{section}{-1}
\begin{abstract}
  We define the notion of a hypercube structure on a functor between two
  strictly commutative Picard categories which generalizes the notion
  of  a cube
  structure on a $G_m$-torsor over an abelian scheme. We use this notion
  to define the intersection bundle\\ $I_{X/S}(L_1,\cdots,L_{n+1})$ 
  of $n+1$ line bundles on a
  relative scheme $X/S$ of relative dimension $n$ and to construct an
  additive 
  structure on the functor $I_{X/S}:PIC(X/S)^{n+1}\F PIC(S)$. Finally,
  we study a 
  section of $I_{X/S}(L_1,\cdots,L_{n+1})$ which  generalizes the
  resultant of $n+1$ polynomials in $n$ variables.
\end{abstract}
\maketitle
\section{Introduction}
\label{a}
\subsection{}
 Soit $G$ un groupe alg{\'e}brique commutatif et $L$ un $G_m$-torseur sur $G$
, le  classique th\'eor\`eme du cube 
(\cite{M3},p58) affirme
que le $G_m$-torseur  suivant sur $G^{3}$ est trivial
\begin{equation}
\label{cub1}  
\t (L) \; = \; m^{\ast }L \wedge
(m_{12}^{\ast }L)^{-1}  \wedge (m_{23}^{\ast }L)^{-1} \wedge (m_{31}^{\ast
}L)^{-1}  \wedge p_{1}^{\ast }L \wedge p_{2}^{\ast }L \wedge p_{3}^{\ast }L 
\end{equation} 
 o\`u $m$ , $m_{ij}$ , $p_{i}$ sont les applications $G\times G
\times G \F G $ d\'efinies par 
 $
m(x_{1},x_{2},x_{3})=x_{i}+x_{j}+x_{k}
$,  
$
m_{ij}(x_{1},x_{2},x_{3})=x_{i}+x_{j}
$ 
 et  
$
p_{i}(x_{1},x_{2},x_{3})=x_{i}
$.

\subsection{} En r{\'e}alit{\'e} si $G$ est une vari{\'e}t{\'e} ab{\'e}lienne, une
trivialisation du 
$G_m$-torseur $\t (L)$ devra v{\'e}rifier des conditions de compatibilit{\'e}
qui s'expriment en termes de  biextensions (notion introduite par
\textsc{Mumford} pour l'{\'e}tude
des groupes formels et {\'e}tudi{\'e}e
par \textsc{Grothendieck} dans un cadre plus g{\'e}n{\'e}ral).
Consid{\'e}rons en effet le $G_m$-torseur $\L =\L (L)$ sur $G \times G$ d{\'e}fini
par:
\[
\L (L) = m^{\ast }L \wedge ( p_1^{\ast }L)^{-1} \wedge ( p_2^{\ast
  }L)^{-1}
\]
Une trivialisation $t$ de $\t (l)$ induit deux lois (afin d'all{\'e}ger les
notations on consid{\'e}rera les fibres de $\L$ au dessus d'un point
g{\'e}n{\'e}rique de $G\times G$):
\[
\begin{array}{lcl}
*_1 : \L _{x,y} \wedge  \L _{x,z}& \F &\L _{x,y+z}\\
*_2 : \L _{x,z} \wedge  \L _{y,z}& \F &\L _{x+y,z}
\end{array}
\]
Ces deux lois v{\'e}rifient alors des propri{\'e}t{\'e}s
d'associativit{\'e}, de 
commutativit{\'e} (que le lecteur pourra {\'e}crire sans
difficult{\'e}) et et de 
compatibilit{\'e} entre elles (qui traduit l'{\'e}galit{\'e} des deux
fa\c{c}ons de d{\'e}velopper $ \L _{x+y,x\pr + y\pr }$ ).
On dit alors que $\L$ est une biextension de $G\times G$ par $G_m$.\\
De plus $\L $ est muni d'un isomorphisme de sym{\'e}trie $s: \tau \ust \L
\simeq \L$ ; on parle alors de biextension sym{\'e}trique.
\subsection{} La notion de structure du cube, introduite par \textsc{Breen}
dans \cite{B2}, 
explicite les notions pr{\'e}c{\'e}dentes. Une struture du cube sur un
$G_m$-torseur $L$ sur $G$ est la donn{\'e}e d'une trivialisation $t$ du 
$G_m$-torseur $\t (L)$ sur $G^3$ telle que les lois partielles $\ast
_1$ et $\ast _2$  induites par $t$ sur $\L (L)$ font de $\L (L)$ une
biextension sym{\'e}trique.
\subsection{} Si $\pi : X \F S $ est une courbe relative lisse  et $J$
d{\'e}signe la composante de degr{\'e} $0$ du sch{\'e}ma de Picard relatif
$\text{PIC} (X/S)$, \textsc{Moret-Bailly} d{\'e}duit de l'existence d'une
structure 
de biextension sur le faisceau de Poincar{\'e} $\mathcal{P}$ sur $J \times
J^{\vee }$ et de l'existence d'une polarisation canonique sur $J$
l'existence d'une biextension $\mathcal{B}$ canonique de $J\times J$
par $G_m$ et il montre ensuite l'existence d'isomorphismes canoniques
\begin{equation} \label{pairing}
\mathcal{B} _{\text{cl} (L) ,\text{cl} (M)} 
\simeq
\drpi (L\tx M)^{-1} \ts \drpi (L)  \ts \drpi (M) \ts  \drpi (\ox )
^{-1}
\end{equation}
\subsection{}\textsc{Deligne} propose dans \cite{d3}  un programme dont une
{\'e}tape est la construction pour tout morphisme projectif plat de
dimension $d$ : $\pi : X \F S$,  du fibr{\'e} d'intersection relativement
{\`a} $S$ , de $d+1$ faisceaux inversibles $L_0 ,
\cdots , L_d$ sur $X$. Il s'agit de construire un $\os$-module inversible   
$ I_{X/S} (L_1 , \cdots , L_{n+1})$ , dont la construction est
fonctorielle en les faisceaux $L_i$ (pour les isomorphismes de
faisceaux) et compatible aux changements de base et qui est
multiplicatif en les faisceaux $L_i$. De fa\c{c}on pr{\'e}cise, on veut
construire un syst{\`e}me d'isomorphismes:
\begin{multline}
 I_{X/S} (L_1 , \cdots , L_i \tx L_i \pr , \cdots  , L_{n+1})
\simeq \\
 I_{X/S} (L_1 , \cdots , L_i , \cdots  , L_{n+1})
\ts
 I_{X/S} (L_1 , \cdots ,  L_i \pr , \cdots  , L_{n+1})
\end{multline}
munis de donn{\'e}es de commutativit{\'e}, d'associativit{\'e} et de
compatibilit{\'e} 
entre ces diff{\'e}rentes lois partielles. Un tel faisceau d'intersection
est construit par \textsc{Deligne} dans \cite{De} dans le cas d'une courbe
relative lisse et cette m{\'e}thode est {\'e}tendue par \textsc{Elkik} dans
 \cite{elkik1} au cas des
morphismes plats de Cohen-Macaulay purement de dimension $n$.\\
\textsc{Moret-Bailly} propose une autre m{\'e}thode dans \cite{MB}
dans le cas des courbes lisses, construisant le
faisceau d'intersection par la formule (\ref{pairing}), bas{\'e}e sur
le fibr{\'e} 
d{\'e}terminant, et montrant ensuite la multiplicativit{\'e} de cette
construction en utilisant la propri{\'e}t{\'e} d'autodualit{\'e} de la
jacobienne. Il explique ensuite comment {\'e}tendre cette construction au
cas de courbes de Cohen-Macaulay.\\
Enfin une autre construction est propos{\'e}e par \textsc{Deligne} dans
\cite{d3},  bas{\'e}e sur des symboles. Cette id{\'e}e est utilis{\'e}e par
\textsc{Aitken} \cite{a}  dans le 
cas d'une courbe singuli{\`e}re quelconque sur une base r{\'e}duite.
\subsection{}
 Dans ce travail on introduit une notion de structure du $p$-cube sur
 un foncteur entre deux cat{\'e}gories de Picard strictement commutatives,
 qui {\'e}tend les d{\'e}finitions de \textsc{Breen}  {\`a} un 
cadre l{\'e}g{\`e}rement plus g{\'e}n{\'e}ral (pour $p=3$, si $G$ est un groupe
 alg{\'e}brique commutatif sur un corps $k$ et $L$ est un $G_m$-torseur
 sur $G$, une structure du 3-cube sur le
 foncteur $\gd :G\F\text{Vect}_k,g\mapsto L_g$ coincide avec la
 d{\'e}finition, donn{\'e}e par \textsc{Breen}, d'une structure du cube sur
 $L$). La donn{\'e}e d'une structure du $p$-cube sur un foncteur $\gd
 :\cc\F\dd$ entre cat{\'e}gories de Picard strictement commutatives munit
la  "dif{\'e}rence sym{\'e}trique $(p-1)$-i{\`e}me
 de $\gd$", qui est un foncteur  $\cc ^{p-1}\F \dd$, de donn{\'e}es
 d'additivit{\'e} en chaque variables.  

On applique ces notions au cas $X/S$ est un sch{\'e}ma relatif projectif
 de dimension $n$ quelconque sur une base localement noeth{\'e}rienne
 et $\gd$ est le foncteur d{\'e}terminant  de l'image directe d{\'e}riv{\'e}e 
$PIC (X/S)\F PIC (S)$. Le th{\'e}or{\`e}me principal de ce travail montre
 alors alors l'existence d'une structure du $(n+2)$-cube canonique sur
 $\gd$. Sa d{\'e}monstration est bas{\'e}e  sur une r{\'e}currence sur
 la dimension  des sch{\'e}mas consid{\'e}r{\'e}s, qui construit une
 structure du $(n+2)$-cube sur le  foncteur d{\'e}terminant
 associ{\'e} {\`a} un sch{\'e}ma relatif de dimension $n$ {\`a} 
partir de la donn{\'e}e de structures du $(n+1)$-cube sur le foncteur
d{\'e}terminant de tout sous-sch{\'e}ma relatif de $X$ de dimension $n-1$. Une
telle r{\'e}currence impose donc, m{\^e}me si on veut montrer le r{\'e}sultat
uniquement pour un
sch{\'e}ma relatif lisse, de savoir le montrer pour des sous-sch{\'e}mas qui
n'ont aucune raison d'{\^e}tre lisses ou m{\^e}me simplement
 r{\'e}duits et donc de travailler avec des sch{\'e}mas $X/S$ assez
 g{\'e}n{\'e}raux. 

Ceci nous permet de montrer que le fibr{\'e}  d'intersection 
$I_{X/S}:PIC (X/S)^{n+1}\F PIC (S)$, d{\'e}fini,  en
 suivant la m{\'e}thode de \textsc{Moret-Bailly}, comme la
diff{\'e}rence sym{\'e}trique $(n+1)$-i{\`e}me du foncteur
d{\'e}terminant est bien muni de donn{\'e}es 
 d'additivit{\'e} en chaque variable. On {\'e}tudie alors une section
 canonique de $I_{X/S}(L_1,\cdots ,L_n)$ qui g{\'e}n{\'e}ralise la notion de
 r{\'e}sultant de $n$ polyn{\^o}mes.

On a fait appel de fa\c{c}on syst{\'e}matique dans les raisonnements {\`a} une
 notion  de $n$-cube dans une cat{\'e}gorie de Picard
strictement commutative. Cette notion est de nature est de nature
combinatoire et permet de repr{\'e}senter de mani{\`e}re commode des
 syst{\`e}mes  d'isomorphismes de la forme $a\otimes b \isom c\otimes
 d$. Son  seul int{\'e}r{\^e}t est de permettre une repr{\'e}sentation 
graphique des raisonnements de r{\'e}currence sur la dimension. Le d{\'e}faut
 de cette approche est qu'il cache l'aspect g{\'e}om{\'e}trique
 li{\'e} {\`a} la  notion de multiextension.
\subsection{Plan de l'article} 
\begin{enumerate}
\item Un exemple des m{\'e}thodes utilis{\'e}es: Les
  propri{\'e}t{\'e}s d'additivit{\'e}  en chaque variable du  "nombre
  d'intersection" de plusieurs diviseurs.
\item Pr{\'e}liminaires techniques sur les diviseurs de Cartier relatifs
  et introduction d'une notion ad hoc de faisceau inversible
  suffisamment positif.
\item Rappels sur les cat{\'e}gories de Picard et introduction de la
  notion de structure du cube.
\item Pr{\'e}sentation des cat{\'e}gories et foncteurs utilis{\'e}s.
\item Th{\'e}or{\`e}me pricipal: l'existence d'une structure du cube sur le
  fibr{\'e} d{\'e}terminant.
\item Applications au fibr{\'e} d'intersection. Construction et {\'e}tude du
  r{\'e}sultant. 
\end{enumerate}
\subsection{Remerciements}
On reconna{\^\i}tra dans ce travail l'influence de Larry
\textsc{Breen}, qui m'a 
introduit dans ce domaine et qui m'a expliqu{\'e} avec patience les
subtilit{\'e}s des structures du cube. Je l'en remercie vivement.
\section{Un exemple introductif}
Soit $X$ un sch{\'e}ma projectif. Pour tout entier $p$ et tous faisceaux
inversibles $L_1 , \cdots , L_p$ sur X, posons:
\[<L_1 , \cdots , L_p>_X =
(-1)^p \chi (\ox ) +
\sum_{k=1}^p (-1)^{p-k}\sum_{1 \leq i_1< \cdots < i_k\leq p}
 \chi (L_{i_1} \otimes \cdots \otimes
 L_{i_k})
\]
Si $p$ est {\'e}gal {\`a} la dimension $n$ de  $X$, on parlera alors de {\em
  nombre d'intersection} de $L_1 , \cdots , L_n$.\\
 Notons d'abord que la d{\'e}finition du nombre d'intersection est
  sym{\'e}trique en les $L_i$ et que pour tout entier $k$ et tous
  faisceaux inversibles  
 $L_1 , \cdots , L_k,L,M$ , on a:
\begin{multline}
<L_1 , \cdots , L_{n-1},L, M>_X =\\
 <L_1 , \cdots , L_{n-1},L>_X
 +<L_1 , \cdots , L_{n-1}, M>_X
 -<L_1 , \cdots , L_{n-1},L \otimes M>_X \; .
\end{multline}
Le r{\'e}sultat suivant est bien connu (cf. par exemple
\cite{beauville},Th 1.4 pour le cas de la dimension 2), mais sa
d{\'e}monstration introduit dans un cadre simple les id{\'e}es
utilis{\'e}es dans ce travail et il sera utilis{\'e}, dans sa reformulation
(\ref{caracteristique2}) pour la d{\'e}monstration du th{\'e}or{\`e}me
principal.  
\begin{lemme} \label{caracteristique}
Soit $X$ un sch{\'e}ma projectif de dimension $n$, alors:
\begin{enumerate}
\item Le morphisme {\em nombre d'intersection}: $PIC(X)^n\F \Z$ est
  $n$-lin{\'e}aire.
\item Si $(\gs _1,\cdots ,\gs _n)$ est une suite r{\'e}guli{\`e}re de sections
  de $L_1 , \cdots , L_n$, alors $<L_1 , \cdots , L_n>_X$ est {\'e}gal {\`a} la
  longueur du sch{\'e}ma $Z$, de dimension 0, des z{\'e}ros de la section 
$\gs =\sum_{i=1}^n \gs _i:\ox\F\bigoplus_{i=1}^{n}L_i$.
\end{enumerate}
\end{lemme}
\begin{proof}[Preuve]
1. Il suffit de montrer par r{\'e}currence sur $n$ 
que si $X$ est un sch{\'e}ma de dimension $n$ et si $L_1 , \cdots , L_{n+1}$
sont  $n+1$ faisceaux inversibles sur
X, on a: $<L_1 , \cdots , L_{n+1}>_X =0$.\\
Cette assertion est {\'e}vidente dans le cas d'un sch{\'e}ma de dimension 0
puisque dans ce cas,
pour tout faisceau inversible $L$, on a $\chi (L)= \text{long}(X)$.\\
Supposons maintenant l'assertion v{\'e}rifi{\'e}e pour tout sch{\'e}ma de
dimension inf{\'e}rieure ou {\'e}gale {\`a} $n$ et consid{\'e}rons un
sch{\'e}ma projectif 
$X$ de dimension $n$. Pour des faisceaux inversibles
$L_1 , \cdots , L_n$ et un diviseur effectif $D$, on a:
\begin{multline}
<L_1 , \cdots , L_n,\ox (D)>_X =(-1)^n (\chi (\ox (D)) - \chi (\ox )\\ +
\sum_{k=1}^n  (-1)^{n-k}\sum_{1 \leq i_1< \cdots < i_k\leq n}
(\chi (L_{i_1} \otimes \cdots \otimes L_{i_k}\otimes \ox (D)) 
-\chi (L_{i_1} \otimes \cdots \otimes L_{i_k}))
\end{multline}
En appliquant la propri{\'e}t{\'e} d'additivit{\'e} de la caract{\'e}ristique
d'Euler-Poincar{\'e}  {\`a} 
des suites exactes de la forme
 $
0 \F L\F L(D)\F L(D)/L \F 0
$
et en identifiant $\chi _X (L(D)/L)$ et  $\chi_D (L(D)|_D)$, on en d{\'e}duit:
\[
<L_1 , \cdots , L_n,\ox (D)>_X =
<L_1(D)|_D , \cdots , L_n(D)|_D>_D
\]
En appliquant l'hypoth{\`e}se de r{\'e}currence au diviseur effectif $D$, on
obtient donc:
\begin{equation}
<L_1 , \cdots , L_n,\ox (D)>_X =0
\end{equation}
On en d{\'e}duit donc,  pour tout diviseur
effectif $D$, les {\'e}galit{\'e}s:
\begin{equation}
\label{add1}
<L_1 , \cdots , L_{n-1},L_n(D)>_X=
<L_1 , \cdots , L_{n-1},L_n>_X 
+
<L_1 , \cdots , L_{n-1},\ox (D)>_X
\end{equation}
et
\begin{multline}
\label{add2}
<L_1 , \cdots ,L_i \otimes L_i\pr ,\cdots , L_{n-1},\ox (D)>_X\\ 
=
<L_1 , \cdots , L_i\pr ,\cdots , L_{n-1},\ox (D)>_X
+
<L_1 , \cdots ,L_i ,\cdots , L_{n-1},\ox (D)>_X
\end{multline}
Comme $X$ est projectif, si $L_n$ est un
faisceau  inversible sur $X$, il peut s'{\'e}crire $\ox (D-E)$, o{\`u}  $D$ et
$E$  sont des diviseurs effectifs. On obtient alors en appliquant (\ref{add1}):
\[
<L_1 , \cdots , L_{n-1},L_n>_X= 
<L_1 , \cdots , L_{n-1},\ox (D)_X>-
<L_1 , \cdots , L_{n-1},\ox (E)>_X
\]
En appliquant (\ref{add2}) {\`a} chacun des termes de droite de
l'{\'e}galit{\'e} pr{\'e}c{\'e}dente, on obtient l'additivit{\'e} de 
 $<L_1 , \cdots ,L_{n-1},L_n>_X$ en chacune des variables $L_1 , \cdots ,
L_{n-1}$, et donc, par sym{\'e}trie, en toutes les variables.\\
2. Consid{\'e}rons le 
complexe de Koszul:
\[
K_{\bullet}:\;
0 \F \L ^n (E) \F \L ^{n-1}(E) \F \cdots \F E \F \ox \F \ox /I_Z \F 0
\]
associ{\'e} au morphisme $\gs\ve : E=\bigoplus_{i=1}^{n}L_i ^{-1} \F
\ox$. Par hypoth{\`e}se, $K_{\bullet}$ est exact, donc $\chi
(K_{\bullet})=0$ et le r{\'e}sultat s'en d{\'e}duit, compte tenu de
l'isomorphisme  
\[
\L ^p (E) \simeq \bigoplus_{1 \leq i_1< \cdots < i_k\leq p}
L_{i_1} \otimes \cdots \otimes
 L_{i_k}.
\]
\end{proof}
Dans la suite de ce travail on {\'e}tudiera, non plus une application
ensembliste de  $Pic(X)^n$ dans $\Z$, mais un foncteur de la cat{\'e}gorie
de Picard $\text{PIC}(X)$ dans une autre cat{\'e}gorie de Picard, et les
difficult{\'e}s proviennent de la n{\'e}cessit{\'e} de faire des constructions
fonctorielles. 
\section{pr{\'e}liminaires techniques}
Dans cette partie nous d{\'e}taillons quelques propri{\'e}t{\'e}s des sections
d'un faisceau inversible sur un sch{\'e}ma relatif et nous introduisons
une notion technique de faisceau suffisamment positif, utile par la
suite. $S$ est ici un sch{\'e}ma localement noeth{\'e}rien et $\pi : X \F S$
est un morphisme projectif et plat.
\subsection{Diviseur relatif d{\'e}fini par une section d'un fibr{\'e}}
\label{div}
Soit $L$ un faisceau inversible sur le sch{\'e}ma relatif $X/S$, une
section  $\s :
\ox \F L$ de $L$ sera dite $\pi$-r{\'e}guli{\`e}re si elle d{\'e}finit elle
d{\'e}finit  un diviseur de Cartier relatif effectif de
$X/S$. $\gs$ est $\pi$-r{\'e}guli{\`e}re si et seulement si pour tout point
$x\in X$, $\s (\oo _{X,x}) 
\subset L_x$ est un $\oo _{X,x}$-module plat et le quotient 
 $(L \tx \oo _{X,x}) / \s (\oo _{X,x})$ est un $\oo_{S,\pi (x)}$-module
 plat.\\
L'ensemble $U\pr$ des points $x$ de $X$ v{\'e}rifiant ces deux
propri{\'e}t{\'e}s 
apparait comme un ensemble de platitude et est donc un ouvert de
$X$ par (\cite{EGA4}, 11.1.1). Notons $Z\pr$ le compl{\'e}mentaire de
$U\pr$. Comme $\pi$ est projectif, $Z= \pi (Z\pr )$ est un ferm{\'e} de $S$
dont nous noterons $U$ le compl{\'e}mentaire.\\
Le th{\'e}or{\`e}me \cite{ma},22.5  entra{\^\i}ne que si $A\F B$ est un
morphisme d'anneaux locaux et $k$ d{\'e}signe le corps r{\'e}siduel de $A$,
 pour tout  $u \in B$, on a l'{\'e}quivalence
des assertions:
\begin{enumerate}
\item $u$ est non diviseur de z{\'e}ro dans $B$ et $B/uB$ est un
  $A$-module plat.
\item $u \otimes _A 1 $ est non diviseur de z{\'e}ro dans $B \otimes _A
  k$.
\end{enumerate}
On d{\'e}duit de ceci qu'un point  $s$ de $S$ est dans $U$ si et seulement
si $\s \mid _{X_s}$ d{\'e}finit un diviseur de Cartier sur $X_s$, soit
encore si $\s $ ne s'annulle sur aucun point associ{\'e} de $X_s$.
\subsection{Suites r{\'e}guli{\`e}res}
Soient $L_1,\cdots , L_p$ des faisceaux inversibles sur  $X$,
muni de sections $\gs_i$.
 Notons $D_i$ le lieu des z{\'e}ros de $\gs _i$. On 
dira que que la suite (ordonn{\'e}e) $(\gs _1,\cdots ,\gs _p)$ est une
{\em suite $\pi$-r{\'e}guli{\`e}re} 
si les deux conditions suivantes sont v{\'e}rifi{\'e}es:
\begin{enumerate}
 \item $D_1\cap\cdots\cap D_p \neq \emptyset$.
 \item Pour tout entier $i$ compris entre 1 et $p$, $\gs _i$ d{\'e}finit
   une section 
   $\pi$-r{\'e}guli{\`e}re sur le sch{\'e}ma relatif
   $D_1\cap\cdots\cap D_{i-1}/S$, 
   si $i>1$, ou sur $X/S$, si $i=1$.
\end{enumerate}
On notera que \begin{enumerate}
\item Si  $(\gs _1,\cdots ,\gs _p)$ est une suite
  $\pi$-r{\'e}guli{\`e}re, alors pour 
  tout $1\leq i \leq p$, $D_1\cap\cdots\cap D_i$ est plat sur $S$.
\item Soit $(\gs _1,\cdots ,\gs _p)$ une suite de sections de $(L_i)$,
 le sous-ensemble de $S$ au dessus duquel  
 $(\gs _1,\cdots ,\gs _p)$ est une suite r{\'e}guli{\`e}re est un
 ouvert de $S$. 
\end{enumerate}

\subsection{Faisceaux inversibles suffisamment positifs}
On dira qu'un faisceau inversible $L$ sur $X$ est 
{\em suffisamment positif} 
(on notera $L \gg 0$) si $L$ est tr{\`e}s ample relativement {\`a} $\pi$ et si
pour tout $i>0$, on a $R^i \pi _{\ast} L =0$. Ces faisceaux v{\'e}rifient
les propri{\'e}t{\'e}s suivantes:
\subsubsection{}  La propri{\'e}t{\'e}  pour un faisceau inversible un
faisceau inversible d'{\^e}tre suffisamment positif est conserv{\'e}e par
 tout changement de base  
 $f:T\F S$. Ceci r{\'e}sulte de l'invariance
  par changement de base de la notion de faisceau relativement tr{\`e}s ample
  (\cite{EGA2}, 4.4.10,iii ) et
  du th{\'e}or{\`e}me de changement de base dans la cohomologie
 (\cite{Ha},12.11). 
\subsubsection{} Si $L_1 , \cdots , \L_k$ sont des
 faisceaux inversibles sur $X$,
  il existe des faisceaux inversibles suffisamment positifs  
 $M_1 , \cdots , M
_k$ tels que les $L_i \tx M_i$ sont tous isomorphes
  et suffisamment positifs. En effet si $M$ est un faisceau inversible
 tr{\`e}s ample relativement {\`a} $\pi$, pour tout faisceau inversible $L$,
 il existe un entier $N$ tel que pour tout $n \geq N$, on a
 $L \tx M^{\otimes n} \gg 0$ (\cite{EGA2},4.4.10.{\em ii} et
  \cite{Ha},th III.8.8.c).  En appliquant ceci
 aux faisceaux 
$\bigotimes_{j\neq i}L_j$ pour $1\leq i\leq k$
 et $ L_1 \otimes \cdots \otimes  L_k $ on trouve
 un entier $n$ tel que les faisceaux 
$\left( \bigotimes_{j\neq i}L_j \right) \otimes M^{\otimes n}$ et
 $\left( \bigotimes_{1\leq j \leq k}L_j \right) \otimes M^{\otimes n}$ soit
 suffisamment positifs. Il suffit alors de prendre 
$M_i =\left( \bigotimes_{j\neq i}L_j
                                \right) \otimes M^{\otimes n}$.
\subsubsection{} \label{cb} Si $L$ est suffisamment positif, $\pi _{\ast} L$
 est un faisceau
 localement libre  sur $S$ et pour tout $s\in S$, la fibre $(\pi _{\ast}
  L)_s$ est isomorphe {\`a} $H^0 (X_s ,L)$  (ceci r{\'e}sulte du
 th{\'e}or{\`e}me de  
changement de base \cite{Ha},12.11). Si de plus $X$ est {\`a} fibres de
 dimension au moins 1,  $\pi _{\ast} L$ est de rang au moins 2.
\subsection{Le diviseur universel d'un faisceau suffisamment positif}
\label{div2}
 Si $L$ est suffisamment positif, $E=(\pi\lst L)\ve$ est un
 $\ox$-module localement libre. Consid{\'e}rons alors le fibr{\'e}
projectif $P_L =  \mathbb{P}(E)$ et effectuons le changement de base:
\[
\begin{CD}
X_{P_L} @>g>> X \\
@V\pi VV         @VV\pi V\\
P_L    @>f>>   S
\end{CD}
\]
Par d{\'e}finition de $P_L$, on a un morphisme surjectif 
$\xymatrix{f\ust E \ar@{->>}[r] &{\oo}_{P_L}(1)}$, 
qui induit un morphisme 
$\xymatrix{{\oo}_{P_L}(-1) \ar@{^{(}->}[r] & (f\ust E)\ve}$. Or on a
$(f\ust E)\ve = f\ust(E\ve)=f\ust(\pi\lst L)=\pi\lst(g\ust L)$, o{\`u} la
derni{\`e}re {\'e}galit{\'e} provient de l'hypoth{\`e}se $L\gg 0$. 
Le morphisme 
$\xymatrix{{\oo}_{P_L}(-1) \ar@{^{(}->}[r] & \pi\lst(g\ust L)}$ 
ainsi obtenu induit par adjonction un morphisme  
$\pi\ust\oo_{P_L}(-1) \F g\ust L$. On obtient donc finalement une section
canonique $\gs_L$ de $g\ust L \otimes \pi\ust\oo_{P_L}(1)$. 
 En appliquant au sch{\'e}ma
relatif $X_{P_L}/P_L$ les constructions de la section (\ref{div}), on
construit un ouvert $U_L$ de $P_L$, au 
dessus duquel $\s _L$ d{\'e}finit un diviseur de Cartier relatif $D_L$ et
on note  $Z_L$  son  compl{\'e}mentaire. On obtient ainsi un isomorphisme
canonique  $(L\otimes\pi\ust\oo_{P_L}(1))| _{X _{U_L}} \simeq
\oo (D_L )$. La situation est d{\'e}crite par le diagramme suivant:
\[ 
\xymatrix{
&(L\otimes\pi\ust\oo_{P_L}(1)\simeq \oo (D_L))\ar@{.}[d] 
&(\ox \stackrel{\gs_L}{\F} L\otimes\pi\ust\oo_{P_L}(1))\ar@{.}[d] 
&L\ar@{.}[d] \\
D_L \ar[dr]  \ar@{^{(}->} [r] &X_{U_L}  \ar[r] \ar[d] &X_{P_L} \ar[r]
\ar[d] 
& X\ar[d]_{\pi} \\
&U_L  \ar@{^{(}->}[r] &P_L \ar[r] &S
}
\]

\begin{lemme} \label{genericite1}
Pour tout point $s$  de $S$, la fibre de $Z_L$ au dessus de $s$ 
est une union finie de sous-espaces
lin{\'e}aires propres  de  $P_{L,s}$.
\end{lemme}
\begin{proof}[Preuve]
En effet, soit $s \in S$, on d{\'e}duit de (\ref{div}) l'{\'e}galit{\'e}:
\[
Z_{L,s} = \bigcup_{x \in \text{Ass}(X_s)}
 \mathbb{P} \left( \left\{\s \in H^0 (X_s ,L) | \s (x)=0 \right\}
 \right)
\subset  \mathbb{P} \left( H^0 (X_s ,L) \right)
= (P_L)_s
\]
L'ensemble des points associ{\'e}s de $X_s$ est fini puisque  $X_s$ est
projectif et comme $L\mid _{X_S}$ est tr{\`e}s ample, pour tout 
$x \in \text{Ass}(X_s)$, l'inclusion
$  \left\{\s \in H^0 (X_s ,L) | \s (x)=0 \right\} \subset  H^0 (X_s
,L) $
est stricte.
\end{proof}
\addtocounter{subsubsection}{1}
\subsubsection{}
Donnons nous de plus une suite $\pi$-r{\'e}guli{\`e}re de sections
$(\gs_i)_{i=1,\cdots ,p}$ de
faisceaux inversibles $L_i$ sur $X$. Soit $V$
l'ouvert  de $U_L$ au dessus duquel  $(\gs_1,\cdots,\gs_p,\gs_L)$ est une
suite $\pi$-r{\'e}guli{\`e}re et notons $Y$ son compl{\'e}mentaire.
\addtocounter{theo}{1}
\begin{lemme}\label{genericite2}
Pour tout point $s$  de $S$, la fibre de $Y$ au dessus de $s$ 
est contenue dans une union finie de sous-espaces
lin{\'e}aires propres  de  $(P_L)_s$.
\end{lemme}
\begin{proof}[Preuve]
Soit $D$ le lieu des z{\'e}ros de la section $\bigoplus\gs_i$ de
$\bigoplus L_i$.
Un point $u\in U_L$ est {\'e}l{\'e}ment de $V$ si et
seulement si $\gs |_{D_u}$ d{\'e}finit un diviseur de Cartier effectif sur
$D_u$. On  raisonne donc comme pour le lemme pr{\'e}c{\'e}dent en {\'e}crivant:
\[
Y_s = \bigcup_{x\in \text{Ass}(D_s)}
\mathbb{P}\{ \gs \in H^0 (X_s ,L) |\gs (x)=0 \}
\]
\end{proof}
\section{Structure du cube sur des  cat{\'e}gories de Picard}
Si $(\mathcal{C},\otimes )$ et  $(\mathcal{D},\otimes ) $ sont des
cat{\'e}gories  de Picard strictement commutatives et $\gd $ est un
foncteur de  $\mathcal{C}$ vers 
$\mathcal{D}$, on d{\'e}finit ici la notion de structure du cube sur $\gd
$.
\subsection{Rappels sur les cat{\'e}gories de Picard commutatives}
\subsubsection{D{\'e}finitions}
Une  {\em cat{\'e}gorie de Picard} est une cat{\'e}gorie $\cc$ dont toutes les
fl{\`e}ches sont des isomorphismes, et qui est  munie d'un bifoncteur
$\otimes : \cc \times \cc \F \cc$ tel que, pour tout $L \in
\text{ob} (\cc )$, les foncteurs 
$X \mapsto L\otimes X $ et $X \mapsto X\otimes L $ soient des
{\'e}quivalences de cat{\'e}gories, et munie de plus de 
donn{\'e}es  d'associativit{\'e} pour
$\otimes$, c'est {\`a} dire d'un syst{\`e}me fonctoriel d'isomorphismes
\[
(L \otimes M) \otimes N \simeq L \otimes (M \otimes N)
\]
v{\'e}rifiant des conditions de compatibilit{\'e} d{\'e}crites par l'axiome du
pentagone.\\
 Une cat{\'e}gorie de Picard $\cc$ est dite {\em commutative} si
 elle est munie 
 de plus de donn{\'e}es de commutativit{\'e}, c'est {\`a} dire d'un syst{\`e}me
 fonctoriel d'isomorphismes 
\[
L\otimes M \isom M \otimes L
\]
compatible avec les donn{\'e}es  d'associativit{\'e} (axiome de
l'hexagone). On prendra garde que l'isomorphisme de commutativit{\'e} 
$L\otimes L \isom L \otimes L$ n'est en g{\'e}n{\'e}ral pas le morphisme
identit{\'e} (quand c'est le cas pour toul $L$, on  dit que la
cat{\'e}gorie est {\em 
  strictement commutative}).\\ 
 On d{\'e}duit des axiomes d'une cat{\'e}gorie de Picard commutative 
l'existence d'un objet unit{\'e} $\oo $, munis de morphismes 
$L \otimes \oo \F L \longleftarrow \oo \otimes L$ compatibles aux
contraintes de commutativit{\'e}.\\
On obtient de m{\^e}me l'existence {\`a} isomorphismes uniques
pr{\`e}s d'objets 
inverses  $L\ve$ munis de morphismes 
$L \otimes L\ve  \F \oo \longleftarrow L\ve \otimes L$.
 \subsubsection{Produit d'une famille index{\'e}e par un ensemble fini}
Soit  $(L_i)_{i\in I}$ une famille d'objets d'une cat{\'e}gorie de Picard
$\cc$ index{\'e}s par un ensemble fini  $I$. Si $I$ est muni d'un ordre
total $<$, les donn{\'e}es 
d'associativit{\'e} de $\cc$ permettent de d{\'e}finir de fa\c{c}on
fonctorielle un objet  
$\bigotimes_{I,<} L_i$ de $\cc$.\\
Supposons maintenant que  $\cc$ est une cat{\'e}gorie de Picard
commutative. Si $<_1$ et $<_2$ sont deux ordres totaux sur $I$, les
donn{\'e}es de commutativit{\'e} de $\cc$ d{\'e}terminent un isomorphisme 
$\bigotimes_{I,<_1} L_i \isom \bigotimes_{I,<_2} L_i$. On peut alors
d{\'e}finir  $\bigotimes_{i\in I}L_i$ comme la limite inductive (ou
projective) des $\bigotimes_{I,<} L_i$ sur tous les ordres totaux $<$
sur $I$.  Pour d{\'e}finir un morphisme 
 $\bigotimes_{i\in I}L_i \F M $ dans $\cc$, il suffira donc de choisir
 un ordre $<$ sur $I$ et de d{\'e}finir un morphisme  
$\bigotimes_{I,<} L_i  \F M$.  Par ailleurs pour deux ensembles
d'indices disjoints $I$ et $J$, on un isomorphisme canonique  
 $\bigotimes_{i\in I}L_i \otimes  \bigotimes_{i\in J}L_i \simeq 
 \bigotimes_{i\in I\cup J}L_i$, obtenu en consid{\'e}rant un ordre sur
 $I\cup J$ tel que $I<J$ et les ordres induits sur $I$ et $J$.
\subsubsection{Foncteur additif.}
Soit $F: \cc \F \dd$ un foncteur additif entre deux cat{\'e}gories de
Picard, une donn{\'e}e d'additivit{\'e} $\mu$ pour $F$ sera la
donn{\'e}e pour tout 
couple d'objets $L,M \in \cc$ d'un isomorphisme fonctoriel en $L$ et
$M$:
\[
\mu _{L,M}: F(L) \otimes F(M) \F  F(L\otimes M)
\]
La donn{\'e}e d'additivit{\'e} $\mu$ sera dite compatible aux donn{\'e}es
d'associativit{\'e} de $\cc$ et $\dd$ si le diagramme suivant, dont les
fl{\`e}ches verticales sont donn{\'e}es par les morphismes
d'associativit{\'e} de 
$\cc$ et $\dd$ 
\[
\begin{CD}
F(L) \otimes (F(M) \otimes F(N)) @>{\text{Id}\otimes\mu}>>F(L) \otimes
F(M \otimes N) 
@>{\mu}>>F(L\otimes (M \otimes N)) \\
@VVV @. @VVV \\
(F(L) \otimes F(M)) \otimes F(N) @>{\mu\otimes\text{Id}}>>
F(L \otimes M) \otimes F(N)
@>{\mu}>> F((L\otimes M) \otimes N)
\end{CD}
\]
est commutatif.\\
Si de plus $\cc$ et $\dd$ sont des cat{\'e}gories de Picard commutatives,
$\mu$ sera dite compatible aux donn{\'e}es de commutativit{\'e} de  $\cc$ et
$\dd$ si le diagramme
\[
\begin{CD}
F(L) \otimes F(M)  @>>> F(L \otimes M)\\
@VVV  @VVV \\
F(M) \otimes F(L)  @>>> F(M \otimes L)
\end{CD}
\]
est commutatif (les
fl{\`e}ches verticales sont donn{\'e}es par les morphismes de
commutativit{\'e} de $\cc$ et $\dd$).
\subsection{n-Cube dans une cat{\'e}gorie de Picard strictement
  commutative} 
On introduit ici une notion de nature combinatoire, qui traduit des
calculs de fractions, du genre $ab ^{-1}= cd ^{-1}$, dans une
cat{\'e}gorie de Picard strictement 
commutative $\cc$, en indexant des objets de $\cc$ par les sommets de
diff{\'e}rents hypercubes de $\R ^n$, ce qui nous permettra dans la suite
de repr{\'e}senter graphiquement certains raisonnements.
\subsubsection{Cubes standard de $\R^n$}
Pour tout entier positif $n$, on consid{\`e}re l'ensemble $C_n=\{0,1\}^n$ des
sommets du n-cube standard de $\R ^n$.\\
 On consid{\'e}rera  pour tout entier $i\leq n$, 
les inclusions de $C_{n-1}$ dans $C_n$:  $\phi _i\pr : (s_1,\cdots
,s_{n-1}) \mapsto (s_1,\cdots ,s_{i-1} ,0,s_{i+1},\cdots ,s_n)$ et 
 $\phi _i\sec : (s_1,\cdots
,s_{n-1}) \mapsto (s_1,\cdots ,s_{i-1} ,1,s_{i+1},\cdots ,s_n)$.
L'image de $C_{n-1}$ par  $\phi _i\pr$ (resp.  $\phi _i\sec$) sera
appel{\'e}e la i-{\`e}me $(n-1)$-face arri{\`e}re (resp. avant) de
$C_n$. De m{\^e}me, 
soient  $\ge_1,\cdots ,\ge_k \in \{ 0,1 \} $ et des indices distincts
$i_1,\cdots,i_k\leq n$, on peut consid{\'e}rer l'inclusion 
$\phi_{i_1=\ge_1,\cdots,i_k=\ge_k}$ de $C_{n-k}$ dans $C_n$.\\
 Pour tout sommet 
$ s=(s_1 , \cdots , s_n)\in C_n$, on notera 
$\ge (s)= (-1)^{n- \sum s_i}$.\\
Enfin il sera commode d'introduire un ordre total sur les sommets de
$C_n$ par:
\[
(s \leq s\pr ) 
\Leftrightarrow
\begin{cases}
\sum s_i < \sum s_i\pr &\\
\text{ou}&\\
(\sum s_i = \sum s_i\pr ) &
 \text{et}\; ( \exists i\leq n ,(s_i > s_i\pr )\; \text{et} \;
(\forall j<i\, , \, s_i=s_i\pr ))
\end{cases}
\]
On notera que l'ordre induit par $\leq $ sur une k-face du cube $C_n$
est encore l'ordre $\leq$ sur $C_k$.
\subsubsection{Arrangements cubiques} On appellera alors {\em
  n-arrangement cubique} dans une cat{\'e}gorie de Picard  
commutative $\cc$
la donn{\'e}e de $2^n$ objets $K_s$ index{\'e}s par les sommets de $C_n$.\\
Toute permutation $\gs\in S_n$ agit sur $\R^n$ par
permutation des coordonn{\'e}es, et induit donc $\gs : C_n \F C_n$. Pour
tout $n$-arrangement cubique $K$ dans $\cc$, on 
 notera alors $\gs\ust K$ le compos{\'e} de $K$ avec la permutation de
 $C^n$ induite par $\gs$.\\
Pour tout  n-arrangement cubique $K$ et tout entier $1\leq i\leq n$,
on consid{\'e}rera alors les $(n-1)$-arrangements cubiques $\fipr_i K$ et
$\fisec_i K$ (i-{\`e}me face arri{\`e}re et avant de $K$). Si $A$ et $B$ sont
deux $(n-1)$-arrangements cubiques, on notera pour tout indice $1\leq
i\leq n$, $(\xymatrix{A \ar@{-}[r]_-i&B})$ le n-arrangement cubique
tel que  $\fipr_i K=A$ et $\fisec_i K=B$. 
De m{\^e}me, si $A_{00}$, $A_{01}$, $A_{10}$ et
$A_{11}$ sont des $(n-2)$-arrangements cubiques, on utilisera, pour deux
indices distincts $1\leq i,j \leq n$, la notation 
$K=
\left(
\begin{array}{c}
\xymatrix{
A_{01}\ar@{-}[d]_-j \ar@{-}[r] &A_{11}\ar@{-}[d] \\
A_{00} \ar@{-}[r]_-i &A_{10}
}
\end{array}
\right)$ pour d{\'e}signer le $n$-arrangement cubique $K$ tel que 
$\phi_{i=\ge ,j=\ge\pr}=A_{\ge ,\ge\pr}, \forall \ge,\ge\pr\in \{0,1\}$.\\
Soient $A$ et $B$ deux $n$-arrangements cubiques tels que
$\fisec_iA=\fipr_iB$. Ecrivons les alors sous la forme:
 $A=
(\xymatrix{U \ar@{-}[r]_-i&V})$ et $B=(\xymatrix{V
   \ar@{-}[r]_-i&W})
$ et posons:
\[
(\xymatrix{U \ar@{-}[r]_-i&V})
\ast_i
(\xymatrix{V\ar@{-}[r]_-i&W})
=(\xymatrix{U\ar@{-}[r]_-i&W})
\]
On dira que $A\ast_iB$ est obtenu par recollement de $A$ et $B$ le
long de leur i-{\`e}me face.
\subsubsection{}
Soit $\gd$ un foncteur de $\cc$ dans une cat{\'e}gorie de Picard
commutative $\dd$, pour tout $n$-arrangement cubique $K$ dans $\cc$,
on pose 
\[
\gd (K) = \bigotimes_{s\in C_n} \gd (K_s)^{\ge (s)}
\]
Si $A$ et$B$ sont des $n$-arrangements cubiques recollables dans la
i-{\`e}me direction, les isomorphismes de commutativit{\'e} et de contraction
dans $\dd$ induisent un isomorphisme canonique
\[
\gd (A\ast_iB)\isom \gd(A) \otimes \gd (B)
\]
 Pour toute permutation $\gs\in S_n$, les isomorphismes de commutativit{\'e} 
dans $\dd$ induisent de m{\^e}me:
\[
\gd (\gs\ust A) \isom \gs (A)
\]
Ces deux isomorphismes sont compatibles entre eux, ce qu'on exprime en
disant que le diagramme suivant est commutatif:
\[
\xymatrix{
\gd (\gs\ust( A\ast_iB))\ar[d]&
 \gd (A\ast_i B)\ar[l] \ar[r] &
\gd (A)\otimes \gd (B)\ar[d]\\
\gd (\gs\ust A\ast_{\gs^{-1}(i)} \gs\ust B)\ar[rr]
&& \gd (\gs\ust A)\otimes\gd (\gs\ust B)
}
\]
\subsubsection{Cubes dans la cat{\'e}gorie de Picard strictement
  commutative $\cc$}
\label{def-cube}
On appellera $1$-cube dans $\cc$ un 1-arrangement cubique quelconque 
$(\xymatrix{L\ar@{-}[r]&M})$.\\
On appellera $2$-cube (carr{\'e}) dans $\cc$ la donn{\'e}e d'un 2-arrangement
cubique $K$ et d'un isomorphisme
 $m_K: \oo \isom \otimes_{s\in C_2} K_s^{\ge (s)}$.\\
On appellera $3$-cube  dans $\cc$ la donn{\'e}e d'un 3-arrangement
cubique 
\[
K=\left(
\begin{array}{c}
\xymatrix{
&K_{001}\ar@{-}[dl]\ar@{-}'[d][dd]\ar@{-}[rr]&
&K_{011}\ar@{-}[dl]\ar@{-}[dd] \\
K_{101}\ar@{-}[dd]\ar@{-}[rr]&&K_{111}\ar@{-}[dd]&\\
&K_{000}\ar@{-}[dl]\ar@{-}'[r][rr]&&K_{010}\ar@{-}[dl]\\
K_{100}\ar@{-}[rr]&&K_{110}
}
\end{array}
\right) ,
\]
 et pour chacune des six 2-faces $F$ de $K$, d'un isomorphisme 
 $m_F: \oo \isom \otimes_{s\in F} K_s^{\ge (s)}$ qui v{\'e}rifient la
 condition de compatibilit{\'e} suivante:
Si $F$ d{\'e}signe l'une des faces de $K$ et $F\pr$ d{\'e}signe la face
oppos{\'e}e, on dispose alors d'un morphisme:
\[
\begin{CD}
\oo @>>> \oo \otimes \oo @>m_F \otimes m_{F\pr }>>
 \bigotimes _{s\in F} K_s^{\ge (s)} 
\otimes 
\bigotimes _{s\in {F\pr}} K_s^{\ge (s)}
@>>>  \bigotimes _{s\in K} K_s^{\ge (s)} 
\end{CD}
\]
On impose que le morphisme ainsi obtenu soit ind{\'e}pendant du choix de la
face $F$.\\
Pour $n\geq 3$, un {\em n-cube}  dans $\cc$  sera la donn{\'e}e d'un
n-arrangement cubique $K$ et d'un isomorphisme $m_F$ associ{\'e} {\`a} chaque
2-face $F$ de $K$, tel que tout sous 3-arrangement cubique de $K$ est un cube.
\subsubsection{Construction de cubes}
(a) Pour tout  $n$-cube  $K$ dans $\cc$, tout objet $L$ de $\cc$ et tout indice
$1\leq i\leq n$, on munit le 
$(n+1)$-arrangement cubique $(\xymatrix{K \ar@{-}[r]_-i&K\otimes
  L})=A$ d'un syst{\`e}me de morphismes $m_F$ associ{\'e}s {\`a}
chaque 2-face de 
$A$, qui en fait un n-cube:\\
Si $F$ n'est pas  parall{\`e}le {\`a} la direction $i$, elle
appartient {\`a} l'un 
des deux $n$-cubes  $K$ ou de $K\otimes L$
et est donc d{\'e}ja muni d'un morphisme $m_F$. Sinon elle s'{\'e}crit
 $F=
\left(
\begin{array}{c}
\xymatrix{Y\ar@{-}[d]_-j\ar@{-}[r] &Y\otimes L\ar@{-}[d]  \\
X\ar@{-}[r]_-i  &X\otimes L
}
\end{array}
\right) $ avec $X,Y\in\text{ob}(\cc )$ et les 
 morphismes de commutativit{\'e} et d'associativit{\'e} de $\cc$ induisent un
 morphisme  
$
m_F:(X\otimes L)\otimes Y\isom X\otimes (Y\otimes L)
$.
On v{\'e}rifie sans peine, en se ramenant au cas o{\`u} $K$ est un
cube, que les relations de compatibilit{\'e} de 
(\ref{def-cube}) entre les $m_F$ sont v{\'e}rifi{\'e}es.\\
(b) Pour tout $n>0$, construisons un  $n$-cube 
$K \isom K_{L_0}(L_1, \cdots ,L_n)$ de la mani{\`e}re suivante:
On pose $K_{L_0}(L_1)= (\xymatrix{L_0 \ar@{-}[r]_-1&L_0\otimes L_1})$
et par r{\'e}currence
 \[
K_{L_0}(L_1, \cdots ,L_{n+1})= 
(\xymatrix{K_{L_0}(L_1, \cdots ,L_n) \ar@{-}[r]_-{n+1}&
 K_{L_0}(L_1, \cdots ,L_n)\otimes L_{n+1}}). 
\] 
 On a alors, pour tout $s \in C_n$ :
\[
( K_{L_0}(L_1, \cdots ,L_n))_s = L_0 \otimes 
\left(
\bigotimes_{i=1}^n L_i^{s_i}
\right) .
\] 
(c) Notons enfin que si $K$ est un $n$-cube, si $K_0$ d{\'e}signe le
$(n-1)$-cube  
$(\phi _i\pr )\ust K$, il existe un objet $L$ de $\cc$ et un
isomorphisme $K \isom (\xymatrix{K_0\ar@{-}[r]_i&K_0\otimes L})$, uniques
{\`a} isomorphisme unique  pr{\`e}s. On en d{\'e}duit donc qu'il
existe des objets 
$L_0, \cdots ,L_n$ et un isomorphisme 
$K \isom K_{L_0}(L_1, \cdots ,L_n)$, uniques {\`a} isomorphisme unique
pr{\`e}s. On dira que l'objet $L_i$ est la {\em i-{\`e}me ar{\^e}te} de K.
\subsubsection{Recollement de cubes}
Consid{\'e}rons deux $n$-cubes $K$ et $K\pr$ dans $\cc$ tels que les deux
sous-cubes $\fisec_i (K)$ et $\fipr_i (K\pr )$ sont {\'e}gaux en tant que 
cubes. Munissons le n-arrangement cubique $K\ast _i K\pr$ d'un syst{\`e}me
de morphismes $m_F$, pour chaque 2-face $F$ de $K\ast _i K\pr$, en faisant un 
 un $n$-cube:
Si $F$ n'est pas parall{\`e}le {\`a} la direction $i$, $F$ est une 2-face de
l'un des deux $n$-cubes $K$ ou $K\pr$ et on prend le $m_F$
correspondant.
Si $F$ est parall{\`e}le {\`a} la direction $i$, elle est obtenue en recollant
une 2-face de $K$ et une de $K\pr$:
\[
\xymatrix{
c\ar@{-}[d] \ar@{-}[r] &d\ar@{}[dr]|{\ast _i}
&d\ar@{-}[d] \ar@{-}[r] &f\ar@{}[dr]|{=}
&c\ar@{-}[d] \ar@{-}[r] \ar@{}[rd]|{F} &f\\
a&b\ar@{-}[u] \ar@{-}[l]_-i &b&e\ar@{-}[u] \ar@{-}[l]_-i 
&a&e\ar@{-}[u] \ar@{-}[l]_-i 
}
\]
et le morphisme $m_F$ est d{\'e}fini par:
\[
\oo \isom \oo \otimes \oo 
\isom 
(a \otimes b\ve  \otimes c\ve  \otimes d )
 \otimes
(b  \otimes e\ve  \otimes d\ve  \otimes f)
\isom  
(a \otimes c\ve  \otimes  c\ve  \otimes f)
\]
Pour montrer que ces $m_F$ v{\'e}rifient les relations de compatibilit{\'e},
il suffit de le faire dans le cas d'un recollement de deux 3-cubes $K$
et $K\pr$ le
long d'une 2-face commune $F$. Cela provient alors imm{\'e}diatement du
fait que, par hypoth{\`e}se, $K$ et $K\pr$ sont des cubes et les morphismes
$m_F$ et $m_F\pr$ associ{\'e}s {\`a} $F$, vu comme face de  $K$ et $K\pr$, sont
les m{\^e}mes.
\subsection{Structure de n-cube}
Soient $\cc$ et $\dd$ des cat{\'e}gories de Picard strictement
commutatives et $\gd :\cc\F\dd$ un foncteur.
\begin{Def}
Une structure du n-cube sur le foncteur $\gd$ est la donn{\'e}e pour tout
n-cube $K$ de $\cc$ d'un morphisme $\psi _K:\oo \isom \gd (K)$ dans
$\cc$ 
v{\'e}rifiant les propri{\'e}t{\'e}s suivantes:
\begin{enumerate}
\item {\em Fonctorialit{\'e}.} Pour tout isomorphisme de n-cubes $f: K
  \isom K\pr$, le diagramme induit 
\[
\xymatrix{
& \gd (K)\ar[dd]^-{\gd (f)} \\
{\oo}  \ar[ur]^-{\psi _K} \ar[dr]_-{\psi _{K\pr}}  \\
& \gd (K\pr )
}
\]
est commutatif.
\item {\em Recollements de cubes.} Soient $K$ et $K\pr$ deux n-cubes
  ayant leur i-{\`e}me (n-1)-face en commun, 
  l'isomorphisme naturel
 $\gd (K) \otimes \gd (K\pr ) \isom \gd (K \ast _i K\pr )$ induit un
  diagramme commutatif:
\[
\xymatrix{
& \gd (K) \otimes \gd (K\pr ) \ar[dd] \\
 {\oo}  \ar[ur]^-{\psi _K \otimes\psi _{K\pr} } 
\ar[dr]_-{\psi _{(K \ast _i K\pr )}}  & \\
&\gd (K \ast _i K\pr )
}
\]
\item {\em Propri{\'e}t{\'e} de sym{\'e}trie.} Pour tout
{\'e}l{\'e}ment $\gs $ de $S_n$, le diagramme suivant, dont la
fl{\`e}che verticale 
est donn{\'e}e par les isomorphismes de commutativit{\'e} de $\dd$,
est commutatif: 
\[
\xymatrix{
& \gd (K) \ar[dd] \\
 {\oo}  \ar[ur]^-{\psi _K} \ar[dr]_-{\psi _{\gs \ust K}} & \\
&\gd (\gs \ust K )
}
\]
\end{enumerate}
\end{Def}
\begin{ex}
 \label{caracteristique2}
Soit $X$ un sch{\'e}ma projectif de dimension $n$, consid{\'e}rons les
cat{\'e}gories $\cc = PIC (X)$ et $\dd = \Z$ (cat{\'e}gorie discr{\`e}te) et le
foncteur $ \gd : \cc \F \dd , L \mapsto \chi (L)$. On a vu (Lemme
\ref{caracteristique}) que
$\gd$ est muni d'une structure du $(n+1)$-cube.
\end{ex}
\subsection{Le n-foncteur multilin{\'e}aire associ{\'e} {\`a} une structure du
  (n+1)-cube}
\subsubsection{}
Consid{\'e}rons un foncteur $\gd : \cc  \F \dd$ entre deux cat{\'e}gories de
Picard strictement commutatives.
 Pour tout  (n-1)-cube $K$ dans $\cc$ et tout
entier  $1\leq
i\leq n$, d{\'e}finissons un foncteur:
\[
\L _{K,i}:  \cc  \F \dd , L \mapsto  \gd
 (\xymatrix{K\ar@{-}[r]_-i&K\otimes L}).
\]
De m{\^e}me, pour tout objet $L$ de $\cc$, d{\'e}finissons un $n$-foncteur:
\[
\L _L :  \cc ^n  \F \dd , (L_1, \cdots , L_n) \mapsto  \gd
(K_L (L_1, \cdots , L_n)).
\]
\addtocounter{theo}{1}
\begin{rem} $\L _L $ est canoniquement muni de donn{\'e}es de sym{\'e}trie 
$\L _L \isom  \gs \ust \L _L$, pour tout $\gs \in S^n$.
\end{rem}
\addtocounter{subsubsection}{1}
\begin{rem}\label{identification}
Pour tout entier $1\leq i\leq n$ on peut identifier canoniquement:
\[
\L _L (L_1, \cdots , L_n) 
\isom 
\L _{K_L(L_1,\cdots ,L_{i-1},L_{i+1},\cdots ,L_n),i} (L_i)
\]
\end{rem}
\addtocounter{subsubsection}{1}
\subsubsection{}
\label{def-multifonct}
\addtocounter{theo}{1}
Donnons nous une structure de $(n+1)$-cube $S$ sur $\gd$. Remarquons d'abord
que $S$ induit, pour tous $L,M \in \text{Ob} (\cc)$,
un isomorphisme canonique de foncteurs $ \L_L \isom \L_M$. En effet,
pour tous $L, L_1, \cdots , L_n \in \text{Ob}(\cc )$, on a
$
K_L (L_1, \cdots , L_n) \simeq K_{\oo} (L_1, \cdots , L_n)\otimes L
$ 
et la structure du cube $S$, appliqu{\'e}e au $(n+1)$-cube 
$
(\xymatrix{
 K_{\oo} (L_1, \cdots , L_n)
\ar@{-}[r]_-{n+1} &
 K_{\oo} (L_1, \cdots , L_n)\otimes L})
$
{\'e}tablit donc un isomorphisme canonique:
\[
\gd ( K_{\oo} (L_1, \cdots , L_n)) 
\isom 
\gd ( K_L (L_1, \cdots , L_n))
\]
On peut donc d{\'e}sormais parler du foncteur $\L$ (en ommettant l'indice
$L$).
\subsubsection{}
\addtocounter{theo}{1}
La donn{\'e}e de la structure du cube sur $\gd$ permet de munir chaque
foncteur  $\L _{K,i}$ d'une donn{\'e}e d'additivit{\'e}:
\[
\mu _{K,i}:\L _{K,i}(L) \otimes \L _{K,i}(M) \isom \L _{K,i}(L\otimes M)
\]
d{\'e}finie par:
\begin{multline} \label{linearite}
\L _{K,i}(L) \otimes \L _{K,i}(M) =
\gd (\xymatrix{K\ar@{-}[r]_-i&K\otimes L})
\otimes
\gd (\xymatrix{K\ar@{-}[r]_-i&K\otimes M}) \isom \\
\gd (\xymatrix{K\ar@{-}[r]_-i&K\otimes L})
\otimes
\gd (\xymatrix{K\ar@{-}[r]_-i&K\otimes M})
\otimes 
\gd \left(
\begin{array}{c}
\xymatrix{K \otimes M \ar@{-}[d]_-{i+1}\ar@{-}[r] 
&K\otimes (L \otimes M) \ar@{-}[d]  \\
K\ar@{-}[r]_-i  &K\otimes L }
\end{array}
\right)\\
\isom \gd (\xymatrix{K\ar@{-}[r]_-i&K\otimes (L\otimes M)})
=\L _{K,i}(L\otimes M)
\end{multline}
En utilisant l'identification \ref{identification}, on voit que $S$
munit $\L$ de $n+1$ donn{\'e}es d'additivit{\'e} partielles:
\begin{multline*}
*_i :\L (L_1, \cdots ,L_{i-1} ,L_i ,L_{i+1}, \cdots ,L_n)
\otimes
\L (L_1, \cdots ,L_{i-1} ,L_i\pr  ,L_{i+1}, \cdots ,L_n)\\
\F
\L (L_1, \cdots ,L_{i-1} ,L_i \otimes L_i\pr
 ,L_{i+1}, \cdots ,L_n).
\end{multline*}
\begin{rem} \label{strict.com}
 On notera que la cat{\'e}gorie de Picard $\dd$ {\'e}tant 
{\em strictement} commutative, le diagramme:
\[
\xymatrix{
(A\ve \otimes A) \otimes A\ve \ar[dr] \ar[dd]\\
& A\ve \\
A\ve \otimes (A \otimes A\ve ) \ar[ur]
}
\]
est commutatif, ce qui permet de ne pas pr{\'e}ciser comment on effectue
les contractions dans (\ref{linearite}).
\end{rem}
\begin{lemme}\label{comm-assoc}
Les donn{\'e}es d'additivit{\'e} $\mu _{K,i}$ pour  $\L _{K,i}$ sont
compatibles aux donn{\'e}es d'as\-so\-cia\-ti\-vi\-t{\'e} et de
commutativit{\'e} de $\cc$ 
et $\dd$.
\end{lemme}
\begin{proof}[Preuve]
Ces deux assertions se traduisent en disant que si $L$, $M$ et $N$ sont des
objets de $\dd$, les isomorphismes canoniques
\[
\gd \left(
\begin{array}{c}
\xymatrix{K \otimes M \ar@{-}[d]_-{i+1}\ar@{-}[r] 
&K\otimes (L \otimes M) \ar@{-}[d]  \\
K\ar@{-}[r]_-i  &K\otimes L }
\end{array}
\right)
\isom
\gd \left(
\begin{array}{c}
\xymatrix{K \otimes L \ar@{-}[d]_-{i+1}\ar@{-}[r] 
&K\otimes (M \otimes L) \ar@{-}[d]  \\
K\ar@{-}[r]_-i  &K\otimes M }
\end{array}
\right)
\]
et
\begin{multline} \label{associativite1}
\gd \left(
\begin{array}{c}
\xymatrix{K \otimes M \ar@{-}[d]_-{i+1}\ar@{-}[r] 
&K\otimes (L \otimes M) \ar@{-}[d]  \\
K\ar@{-}[r]_-i  &K\otimes L }
\end{array}
\right)
\otimes
\gd \left(
\begin{array}{c}
\xymatrix{K \otimes N \ar@{-}[d]_-{i+1}\ar@{-}[r] 
&K\otimes (L \otimes M\otimes N) \ar@{-}[d]  \\
K\ar@{-}[r]_-i  &K\otimes (L\otimes M) }
\end{array}
\right)\\
\isom
\gd \left(
\begin{array}{c}
\xymatrix{K \otimes N \ar@{-}[d]_-{i+1}\ar@{-}[r] 
&K\otimes (M \otimes N) \ar@{-}[d]  \\
K\ar@{-}[r]_-i  &K\otimes M }
\end{array}
\right)
\otimes
\gd \left(
\begin{array}{c}
\xymatrix{K \otimes (M\otimes N) \ar@{-}[d]_-{i+1}\ar@{-}[r] 
&K\otimes (L \otimes M\otimes N) \ar@{-}[d]  \\
K\ar@{-}[r]_-i  &K\otimes L }
\end{array}
\right)
\end{multline}
identifient les trivialisations de chacun des termes, donn{\'e}es par la
structure du (n+1)-cube.\\
Le premier point provient de l'hypoth{\`e}se que la structure du cube est
sym{\'e}trique et de la remarque (\ref{strict.com}).\\
Pour l'associativit{\'e}, remarquons d'abord que chacun des membres de
\ref{associativite1}  est isomorphe {\`a}
\begin{multline} \label{associativite3}
\gd \left(
\begin{array}{c}
\xymatrix{K \otimes N \ar@{-}[d]_-{i+1}\ar@{-}[r] 
&K\otimes (M \otimes N) \ar@{-}[d]  \\
K\ar@{-}[r]_-i  &K\otimes M }
\end{array}
\right)
\otimes
\gd \left(
\begin{array}{c}
\xymatrix{K \otimes (M\otimes N) \ar@{-}[d]_-{i+1}\ar@{-}[r] 
&K\otimes (L \otimes M\otimes N) \ar@{-}[d]  \\
K\otimes M \ar@{-}[r]_-i  &K\otimes (L\otimes M) }
\end{array}
\right)\\
\otimes
\gd \left(
\begin{array}{c}
\xymatrix{K \otimes M \ar@{-}[d]_-{i+1}\ar@{-}[r] 
&K\otimes (L \otimes M) \ar@{-}[d]  \\
K\ar@{-}[r]_-i  &K\otimes L }
\end{array}
\right).
\end{multline}
Consid{\'e}rons alors l'empilement de (n+1)-cubes suivant:
\begin{equation}
\label{associativite2}
\xymatrix{
K \otimes N \ar@{-}[d]_-{i+1}\ar@{-}[r] 
&K\otimes (M \otimes N) \ar@{-}[d]_-{i+1} \ar@{-}[r] 
&K\otimes (L\otimes M \otimes N) \ar@{-}[d]   \\
K \ar@{-}[r]_-i 
&K\otimes M \ar@{-}[d] \ar@{-}[r]_-i \ar@{-}[d]_-{i+1} 
&K\otimes (L\otimes M ) \ar@{-}[d]   \\
&K \ar@{-}[r]_-i &K\otimes L
}
\end{equation}
Il r{\'e}sulte alors de la  propri{\'e}t{\'e} 2. d'une structure du cube,
appliqu{\'e}e aux deux fa\c{c}ons d'effectuer des recollements dans
(\ref{associativite2})  que  les
trivialisations 
des membres de (\ref{associativite1}) sont identifi{\'e}es  aux
trivialisations de (\ref{associativite3}).
\end{proof}
De ces diff{\'e}rents r{\'e}sultats, on d{\'e}duit la:
\begin{prop}
La donn{\'e}e d'une structure du (n+1)-cube $S$ sur $\gd: \cc \F \dd$
munit le $n$-foncteur associ{\'e} $\L : \cc ^n \F \dd$ de donn{\'e}es
d'additivit{\'e} $\ast _i$ en chacune des $n$ variables. Ces donn{\'e}es
sont compatibles aux donn{\'e}es d'associativit{\'e}  et de
commutativit{\'e} de 
$\cc$ et $\dd$ ainsi qu'aux donn{\'e}es de sym{\'e}trie
de $\L$ et sont compatibles entre elles.
\end{prop}
\begin{proof}[Preuve]
Les questions d'associativit{\'e} et de commutativit{\'e} proviennent des
r{\'e}sultats analogues pour $\L _{K,i}$ (lemme \ref{comm-assoc}).\\
Pour simplifier les notations exprimons la compatibilit{\'e} des $\ast _i$
entre elles dans le cas $n=2$. Cela se traduit par la commutativit{\'e} du
diagramme:
\[
\begin{CD}
\L (L,N) \otimes \L (L,P) \otimes \L (M,N) \otimes \L (M,P)
@>{\ast _1 \otimes \ast _1}>>
\L (LM,N) \otimes \L (LM,P)\\
@V{\ast _2 \otimes \ast _2}VV @V{\ast _2 }VV \\
\L (L,NP) \otimes \L (M,NP) 
@>{\ast _1}>>
\L (LM,NP)
\end{CD}
\]
\textsc{Breen} montre dans  \cite{B2},2.5 que la commutativit{\'e} de ce
diagramme est 
une cons{\'e}quence de l'associativit{\'e} de $\ast _1$ et $\ast
_2$. En effet, 
cela se traduit en disant que l'isomorphisme canonique d{\'e}duit des
morphismes de contraction:
\begin{multline*}
\ga : \gd (K(LM,N,P)) \otimes \gd (K(L,M,N)) \otimes \gd (K(L,M,P))\\
\isom
\gd (K(L,M,NP)) \otimes \gd (K(L,N,P)) \otimes \gd (K(M,N,P))
\end{multline*}
identifie les trivialisations de chacun des deux termes qui sont
d{\'e}duites de la structure du cube. Cette assertion provient alors du
fait que $\ga$ se d{\'e}compose en
\[
\begin{CD}
\gd (K(LM,N,P)) \otimes \gd (K(L,M,N)) \otimes \gd (K(L,M,P))\\
@VVV\\
\gd (K(L,MN,P)) \otimes \gd (K(M,N,P)) \otimes \gd (K(L,M,N))\\
@AAA\\
\gd (K(L,M,NP)) \otimes \gd (K(L,N,P)) \otimes \gd (K(M,N,P))
\end{CD}
\]
provenant des morphismes d'associativit{\'e} de $\ast _1$ et $\ast _2$.\\
Pour traduire la compatibilit{\'e} des donn{\'e}es d'additivit{\'e}  avec les
donn{\'e}es de sym{\'e}trie, consid{\'e}rons des objets $L_1, \cdots ,
L_n,L,L\pr$ de 
$\cc$ et $\gs \in S_n$ et notons 
$\underline{L}_i= (L_1,\cdots,L_{i-1},L,L_{i+1},\cdots ,L_n)$,
$\underline{L\pr}_i= (L_1,\cdots,L_{i-1},L\pr,L_{i+1},\cdots ,L_n)$ et 
$\underline{L\sec}_i= (L_1,\cdots,L_{i-1},L\otimes L\pr
,L_{i+1},\cdots ,L_n)$ et $j=\gs ^{-1} (i)$.
 On veut montrer que le diagramme suivant, dont les fl{\`e}ches verticales
 proviennent des donn{\'e}es de sym{\'e}trie de $\L$
\[
\begin{CD}
\L (\gs \ust \underline{L}_i )
\otimes 
\L (\gs \ust \underline{L\pr}_i ) @>{\ast _j}>>
\L (\gs \ust \underline{L\sec}_i ) \\
@VVV @VVV \\
\L ( \underline{L}_i )
\otimes 
\L ( \underline{L\pr}_i ) @>{\ast _j}>>
\L ( \underline{L\sec}_i )
\end{CD}
\]
est commutatif. Ceci provient de la propri{\'e}t{\'e} de sym{\'e}trie de la
structure de $(n+1)$-cube, appliqu{\'e}e au  $(n+1)$-cube 
$C= (K(L_1,\cdots,L_{i-1},L, L\pr,L_{i+1},\cdots ,L_n)$ et {\`a} la
permutation $\tau \in S_{n+1} $ telle que 
$\tau \ust C =
K(L_{\gs (1)},\cdots,L_{\gs (i-1)},L, L\pr,L_{\gs (i+1)},\cdots
,L_{\gs (n)})$.
\end{proof}
\section{Les cat{\'e}gories de Picard consid{\'e}r{\'e}es}
Dans ce chapitre, on consid{\`e}re un morphisme projectif et plat $\pi :
X \F S$   sur un sch{\'e}ma localement noeth{\'e}rien et on introduit les
cat{\'e}gories de Picard qui nous serviront dans la suite.
\subsection{}
 $\cc$ d{\'e}signera la
cat{\'e}gorie $PIC (X)$ dont les objets sonts les faisceaux inversibles
sur $X$, les fl{\`e}ches sont les isomorphismes de $\ox$-modules,
$\otimes$ d{\'e}signe le produit tensoriel usuel de deux $\ox$-modules et
les morphismes d'associativit{\'e} et de commutativit{\'e} sont ceux
usuels. Les objets unit{\'e}s et inverses seront simplement le faisceau
structural $\ox$ et le faisceau dual $L^{-1}$.
\subsection{}
 $\dd$ d{\'e}signera la
cat{\'e}gorie $PICgr (S)$ dont les objets sont les couples $(L,d)$ form{\'e}s
d'un faisceau inversible $L$ 
sur $S$ et d'une application localement constante $d: X \F \Z$, les
fl{\`e}ches entre deux objets $(L,d)$ et $(M,e)$ n'existent que si $d=e$
et sont dans ce cas  les isomorphismes entre les  $\os$-modules $L$ et
$M$. Le produit tensoriel sera alors d{\'e}fini par 
\[
(L,d) \otimes (M,e) = (L\tx M , d+e).
\]
Les morphismes de commutativit{\'e} 
$\psi : (L,d) \otimes (M,e) \F  (M,e) \otimes (L,d)$ 
sont donn{\'e}s par:
$\psi : L \ts M \F M \ts L : l \otimes m \mapsto (-1)^{d.e} m\otimes
l$. On notera que la cat{\'e}gorie  $PICgr (S)$ n'est pas strictement
commutative.\\
L'objet unit{\'e} de $\dd$ sera alors $\oo = (\os ,0)$ et l'inverse de
$(L,d)$ sera $(L\inv ,-d)$. Le morphisme d'{\'e}valuation 
$(L,d) \otimes (L\inv ,-d)\F \oo$ est donn{\'e} par l'{\'e}valuation usuelle 
$L\ts L\inv  \simeq \os$. On prendra garde 
 que le morphisme 
$(L\inv ,-d) \otimes (L,d)\F \oo$ est donn{\'e} par le morphisme usuel
d'{\'e}valuation multipli{\'e} par $(-1)^d$.\\
On consid{\`e}rera enfin le foncteur $\gd : \cc \F \dd$ donn{\'e} par le
d{\'e}terminant de l'image  directe d{\'e}riv{\'e}e dont l'existence,
annonc{\'e}e par \textsc{Grothendieck}, est
montr{\'e}e par \textsc{Knudsen} et \textsc{Mumford} dans \cite{KM}.
 $\gd$  associe {\`a} tout $\ox$-module inversible
$L$,  le $\os$-module $\drpi L $, gradu{\'e}
par la fonction localement constante $ s \mapsto \chi (L |
_{X_s})$.\\
En r{\'e}alit{\'e}, le foncteur $\gd$ est d{\'e}fini
sur la cat{\'e}gorie plus {\'e}tendue $COH (X/S)$ des $\ox$-modules
coh{\'e}rents 
et plats sur $S$. Il v{\'e}rifie les propri{\'e}t{\'e}s suivantes:
\subsubsection{}
Toute suite exacte dans  $COH (X/S)$:
\[
0 \F E \F F \F G \F 0
\]
induit un isomorphisme canonique de multiplicativit{\'e}  
$\gd (E) \otimes \gd (G) \isom \gd (F)$
 et pour tout  diagramme de suites exactes courtes:
\[
\begin{CD}
@. 0 @. 0 @. 0 \\
 @.@VVV   @VVV      @VVV\\
0 @>>> E\pr @>>>F\pr @>>> G\pr @>>>0\\
@. @VVV   @VVV      @VVV\\
0 @>>> E @>>>F @>>> G @>>>0\\
@. @VVV   @VVV      @VVV\\
0 @>>> E\sec @>>>F\sec @>>> G\sec @>>>0\\
@. @VVV   @VVV      @VVV\\
@. 0 @. 0 @.0 
\end{CD}
\]
le diagramme
\begin{equation}
\label{diagramme-des-neuf}
\begin{CD}
\gd (E\pr ) \otimes \gd (E\sec ) \otimes \gd (G\pr )\otimes \gd (G\sec
 )
 @>>> \gd (E) \otimes \gd (G)\\
@VVV   @VVV      \\
 \gd (F\pr ) \otimes \gd (F\sec ) @>>> \gd (F)
\end{CD}
\end{equation}
qu'on en d{\'e}duit par application des morphismes de multiplicativit{\'e} et
de commutativit{\'e}  est commutatif
(\cite{KM},prop.1).
Notons que c'est la d{\'e}finition des isomorphismes de commutativit{\'e} dans
$\dd$ qui rend possible la commutativit{\'e} du diagramme pr{\'e}c{\'e}dent.
\subsubsection{}Si $E$ est un $\ox$-module coh{\'e}rent et plat sur $S$, {\`a}
support dans un sous-sch{\'e}ma $Y$ de $X$, on a: $\drpi (E) = \det
\text{R} (\pi _{Y/S}) \lst  (E|_Y)$. 

\subsubsection{} La formation de $\drpi$ commute aux changements de base.
\subsection{} On consid{\`e}rera enfin la cat{\'e}gorie de Picard
strictement 
commutative $\ddp =PIC (S)$ et on notera  $\gdp$ le compos{\'e} du
foncteur $\gd = \drpi : \cc \F \dd$ avec le foncteur  oubli de la
graduation $\dd \F \ddp$. On notera qu'en travaillant dans cette
cat{\'e}gorie $\ddp$ on gagne le fait qu'elle est strictement
commutative, mais 
on perd la possibilit{\'e} d'{\'e}crire certains diagrammes
commutatifs du 
paragraphe pr{\'e}c{\'e}dent, comme par exemple (\ref{diagramme-des-neuf}),
qui ne s'expriment naturellement que gr{\^a}ce au foncteur $\gd$.
\subsection{Structure du cube dans la cat{\'e}gorie des faisceaux
  inversibles gradu{\'e}s} 
Pour montrer l'existence {\'e}ventuelle d'un structure du $p$-cube sur le
foncteur $\gdp : \cc \F \ddp$ entre
cat{\'e}gories de Picard strictement commutatives, il sera n{\'e}cessaire de
passer par l'interm{\'e}diaire du foncteur  $\gd : \cc \F \dd$ {\`a}
valeurs 
dans une cat{\'e}gorie de Picard non strictement commutative. Examinons
ici comment les axiomes d'une structure de $p$-cube sur  $\gdp$ se
traduisent en termes de $\gd$.
\begin{nota}
Si $K$ est un $p$-cube dans $\cc$, pour tout couple d'indices
distincts $i$ et $j$, les isomorphismes de commutativit{\'e} dans $\dd$
induisent un isomorphisme $\gd (\gs _{ij}\ust K)\F \gd (K)$ dans
$\dd$, qui induit, par oubli de la graduation, un isomorphisme 
 $\gdp (\gs _{ij}\ust K)\F \gdp (K)$ dans
$\ddp$. Celui ci diff{\`e}re de celui induit par les isomorphismes de
commutativit{\'e} dans $\ddp$ par un signe, que l'on notera $\ge _{ij}(K)$. 
\end{nota}
\begin{rem}
Si $X\F S$ est {\`a} fibres de dimension $n$, pour tout $(n+2)$-cube 
\begin{equation}
\label{ecriture-cube}
K=\left(
\begin{array}{c}
\xymatrix{C\ar@{-}[d]_-j\ar@{-}[r] 
&D \ar@{-}[d]  \\
A\ar@{-}[r]_-i  &B }
\end{array}
\right) ,
\end{equation}
 on a 
$\chi (A) = \chi (B) =
\chi (C) = \chi (D)$ 
 (d'apr{\`e}s l'exemple
\ref{caracteristique2}) et de plus 
 $\ge _{ij}(K)=(-1)^{\chi (A)}$.
\end{rem}
\addtocounter{subsubsection}{2}
\subsubsection{}
\addtocounter{theo}{1}
Soit $K$ un $p$-cube dans $\cc$ et soit  $t: \os \isom \gdp (K)$ une
 trivialisation de $\gdp (K)$. Pour tout entier $1\leq i\leq p$, $t$
 induit un isomorphisme
$t_i: \gdp (\fipr_i K) \isom \gdp (\fisec_i K)$ dans
 $\ddp$. 
 Le choix de l'ordre usuel sur les sommets de
$\fipr_i K$ et $\fisec_i K$
  associe {\`a}   $t_i$ un
 isomorphisme 
$s_i: \gd (\fipr_i K) \isom \gd (\fisec_i K)$
 dans $\dd$. Exprimons alors, {\`a} l'int{\'e}rieur de la
 cat{\'e}gorie $\ddp$ les 
 relations entre $s_i$ et $s_j$, pour $i\neq j$. On peut associer
 {\`a} $s_i$ le 
 morphisme dans $\dd$:
\[
\begin{CD}
\overline{s_i}:\oo @>>>
(\gd  (\fipr_i K))^{-1} \otimes \gd ( \fipr_i K) 
@>{\text{Id}\otimes s_i}>> 
(\gd  (\fipr_i K))^{-1} \otimes \gd ( \fisec_i K)
@>>>
\gd (K)
\end{CD}
\]
\begin{lemme}
La donn{\'e}e d'une trivialisation $t$ de $\gdp(K)$ est 
{\'e}quivalente {\`a} la donn{\'e}e d'une collection $(s_i)_{1\leq i\leq p}$ 
d'isomorphismes 
$s_i: \gd (\fipr_i K) \isom \gd (\fisec_i K)$ dans $\dd$ telle que, pour
deux indices $i$ et $j$ distincts, les trivialisations 
$\overline{s_i}$ et $\overline{s_j}$ 
de $\gd (K)$ induites par $s_i$ et $s_j$
diff{\`e}rent d'un signe $\ge _{ij}(K)$.
\end{lemme}
\addtocounter{subsubsection}{1}
\begin{proof}[Preuve]
Reprenons la notation (\ref{ecriture-cube}). Le 
 diagramme, dont les fl{\`e}ches verticales sont
donn{\'e}es par les isomorphismes les isomorphismes structuraux de  $\ddp$
 est commutatif:
\[
\begin{CD}
(\gdp (A) \ts \gdp (C)\ve ) \ts (\gdp (C)\ts\gdp (D))
@>{t_i\otimes \text{id} }>>
(\gdp (B) \ts \gdp (D)\ve ) \ts (\gdp (C)\ts\gdp (D))\\
@VVV @VVV \\
(\gdp (A) \ts \gdp (D) ) \ts (\gdp (C)\ts\gdp (C)\ve )
@.
(\gdp (B) \ts \gdp (C)\ve ) \ts (\gdp (D)\ts\gdp (D)\ve )\\
@VVV @VVV \\
(\gdp (A) \ts \gdp (D) )  @>t>> (\gdp (B) \ts \gdp (C) ) \\
@AAA @AAA \\
(\gdp (A) \ts \gdp (D) ) \ts (\gdp (B)\ts\gdp (B)\ve )
@.
(\gdp (C) \ts \gdp (B)\ve ) \ts (\gdp (D)\ts\gdp (D)\ve )\\
@AAA @AAA \\
(\gdp (A) \ts \gdp (B)\ve ) \ts (\gdp (B)\ts\gdp (D))
@>{t_j\otimes \text{id} }>>
(\gdp (C) \ts \gdp (D)\ve ) \ts (\gdp (B)\ts\gdp (D))
\end{CD}
\]
Consid{\'e}rons le diagramme obtenu en rempla\c{c}ant dans le
pr{\'e}c{\'e}dent $\gdp$ par 
$\gd$ et en utilisant les morphismes structuraux de $\ddp$. Son circuit
ext{\'e}rieur est donc 
commutatif {\`a} un signe $(\ge _{ij})^k$ pr{\`e}s, o{\`u} $k$ est le nombre de
transpositions apparaissant dans le 
diagramme. On constate alors que ce nombre est impair.\\
R{\'e}ciproquement, si $(s_i)_{1\leq i\leq p}$ est une telle collection de
morphismes, chaque $s_i$ induit un isomorphisme 
$t_i:\gdp \left(\fipr_i K \right)
\isom 
\gdp\left( \fisec_i K \right)$ dans $\ddp$. La trivialisation de $\gdp
(K)$ dans $\ddp$ est alors ind{\'e}pendante de $i$.
\end{proof}
\subsubsection{Recollements de cubes}
\addtocounter{theo}{1}
Soient  $K$  et $K\pr$ deux $p$-cubes dans $\cc$ et $i$ un indice tel que
$\fipr_i K = \fisec_i K$. Si   $t , t\pr$ sont des
trivialisations de $\gdp (K)$ et $\gdp (K\pr )$, elles induisent une
trivialisation  
$t \ast_i t\pr: \ox \isom \ox \otimes \ox 
\isom \gdp (K) \otimes \gdp (K\pr ) \isom \gdp (K\ast_i K\pr )$.
Notons 
$s_i: \gd (\fipr_i K) \isom \gd (\fisec_i K)$, 
$s\pr_i: \gd (\fipr_i K\pr ) \isom \gd (\fisec_i K\pr )$ et 
$s\sec_i: \gd (\fipr_i (K\ast _i K\pr )) 
\isom \gd (\fisec_i
(K\ast _i K\pr ))$ les isomorphismes dans $\ddp$ associ{\'e}s {\`a} $t$, $t\pr$
et $t \ast_i t\pr$. On a alors $s\sec_i =s\pr_i \circ s_i$. 

On peut alors regrouper les r{\'e}sultats de ce paragraphe dans la
\begin{prop} \label{cubebis}
Une structure de $p$-cube sur $\gdp$ est {\'e}quivalente {\`a} la donn{\'e}e,
pour tout $p$-cube $K$ dans $\cc$, de $p$ isomorphismes dans $\dd$:
\[
s_{K,i}:\gd \left(\fipr_i K \right)
\isom 
\gd\left( \fisec_i K \right)
\]
tels que:
\begin{enumerate}
\item  Pour tout isomorphisme de $p$-cubes $f: K \isom K\pr$, le
  diagramme induit 
\[
\begin{CD}
\gd \left( \fipr_i K \right)
@>>> 
\gd\left( \fisec_i K \right)\\
@VVV @VVV\\
\gd \left( \fipr_i K\pr \right)
@>>> 
\gd\left( \fisec_i K\pr \right)
\end{CD}
\]
est commutatif.
\item Soient $K$ et $K\pr$ deux $p$-cubes ayant leur i-{\`e}me $(p-1)$-face
  en commun, on a l'{\'e}galit{\'e}:
\[
s_{K\ast _i K\pr ,i} = s_{ K\pr ,i} \circ s_{K ,i}
\]
\item Les trivialisations de $\gd (K)$ induites par $s_{K,i}$ et $s_{K,j}$
diff{\`e}rent d'un signe $\ge _{ij}(K)$.
\item Pour toute permutation $\gs _{ij} \in S_p$ de deux indices
  distincts $i$ et $j$, on a:
\[
s_{\gs _{ij}\ust (K),i} = s_{K,i}
\]
o{\`u}  d{\'e}signe la permutation des indices $i$ et $j$.
\end{enumerate}
\end{prop}
\section{Constructions de structures du cube}
\subsection{Cas de la dimension 0} \label{norme}
Si $\pi : X \F S$ est un morphisme fini et plat, le foncteur $\gd =\drpi$ se
r{\'e}duit au foncteur:
\[
\gd : PIC (X) \F PICgr (S)\; ,\;  L
 \mapsto 
\left(\det (\pi \lst L), \deg \pi  \right)
.
\]
et $\gd \pr : PIC (X) \F PIC (S)$ est simplement le foncteur: 
$ L\mapsto \det (\pi \lst L)$.
\subsubsection{Norme }
Rappellons quelques propri{\'e}t{\'e}s de la norme  (cf \cite{EGA2},6.5 et 
\cite{FD1},3.1) pour
un morphisme fini et plat. 
Notons d'abord que, comme $X$ est fini sur $S$, pour tout faisceau
inversible $L$ sur $X$ et  tout $s\in S$, il
existe un ouvert 
$U$ de $S$ contenant $s$ tel que $L$ est trivial sur $X_U$. Si $\ga$
est une section inversible de $\ox$, c'est {\`a} dire un automorphisme de
$\ox$, la norme $N_{X/S} (\ga )$ est le d{\'e}terminant de l'automorphisme
de $\det \pi \lst \ox$ induit par $\ga$. Enfin la norme du faisceau
inversible $L$ est par d{\'e}finition le $\os$-module inversible
\[
N_{X/S} (L) =
(\det \pi \lst L) \ts (\det \pi \lst \ox )^{-1}
\]
Si $(U_i)_{i\in I}$ est un recouvrement ouvert de $S$ tel que $L$ est
trivial sur $X_{U_i}$ avec des fonctions de transition $(g_{ij})_{i,j
  \in I}$, alors
$N_{X/S} (L)$ est trivial sur chaque $U_i$ et a pour fonctions de
transition les $N_{X/S}(g_{ij})$.
\subsubsection{}
Construisons une structure du 2-cube (carr{\'e}), qui traduit les
propri{\'e}t{\'e}s de la norme pour un morphisme fini et
plat:
Il s'agit de construire, pour tous faisceaux inversiblest $L$, $M$,
$N$ et $P$  sur $X$ et  tout 
isomorphisme $\phi : L \tx M \isom N \tx P$, un isomorphisme 
\[
(\det \pi \lst L) \ts (\det \pi \lst M) 
\isom
 (\det \pi \lst N) \ts (\det \pi \lst P) .
\]
D'apr{\`e}s les remarques de la
 section pr{\'e}c{\'e}dente, il suffit de construire pour tout isomorphisme 
$\phi : \ox \tx \ox  \isom \ox \tx \ox$ un isomorphisme 
$ \gd _{\phi}:\det (\pi \lst  \ox ) \ts \det (\pi \lst  \ox )
 \isom 
\det (\pi \lst  \ox ) \ts \det (\pi \lst  \ox )$ tel que si $\phi$ et
 $\psi$ sont deux tels isomorphismes et 

 $\ga$, $\gb$, $\gc$ et $\gd$ des automorphismes de $\ox$ tels que le
 diagramme 
\[
\begin{CD}
 \ox \tx \ox  @>\phi >> \ox \tx \ox \\
@V{\ga \otimes \gb}VV @VV{\gc \otimes \gd}V\\
 \ox \tx \ox  @>\psi >> \ox \tx \ox 
\end{CD}
\]
est commutatif, alors le diagramme induit:
\[
\begin{CD}
 (\det \pi \lst \ox ) \ts (\det \pi \lst \ox ) 
 @>{\det \pi \lst\phi}>> 
(\det \pi \lst \ox ) \ts (\det \pi \lst \ox ) \\
@V{N_{X/S} (\ga ) \otimes N_{X/S}(\gb )}VV 
@VV{N_{X/S}(\gc ) \otimes N_{X/S}(\gd )}V \\
(\det \pi \lst \ox ) \ts (\det \pi \lst \ox ) 
 @>{\det \pi \lst\psi}>> 
(\det \pi \lst \ox ) \ts (\det \pi \lst \ox )
\end{CD}
\]
est commutatif. 
$\phi$ est un automorphisme de $ \ox \tx \ox$, c'est {\`a} dire une
section inversible de $\ox$. On peut consid{\'e}rer la section inversible 
$N_{X/S}(\phi )$ de $\os$, qui d{\'e}finit donc  l'automorphisme recherch{\'e}
de  $(\det \pi \lst \ox ) \ts (\det \pi \lst \ox )$. La condition de
fonctorialit{\'e} se traduit en disant que, si $\ga , \gb , \gc ,\gd , \phi
,\psi$ sont des sections inversibles de $\ox$ telles que 
$\ga\gb\psi =\gc\gd\phi$, alors 
\[
N_{X/S}(\ga ) N_{X/S}(\gb) N_{X/S}(\psi ) 
 = N_{X/S}(\gc ) N_{X/S}(\gd ) N_{X/S}(\phi )\; ,
\]
 ce qui est simplement la
 mutiplicativit{\'e} de la norme.\\
Le lemme suivant est une cons{\'e}quence imm{\'e}diate de la
 d{\'e}finition de la 
 norme d'un faisceau inversible et des propri{\'e}t{\'e}s de
 multiplicativit{\'e} 
 de la norme.
\addtocounter{theo}{2}
\begin{lemme}
\begin{enumerate}
\item La construction pr{\'e}c{\'e}dente d{\'e}termine une structure
  du carr{\'e}  sym{\'e}trique sur $\gdp$.
\item Le foncteur lin{\'e}aire $PIC(X) \F PIC (S)$ d{\'e}duit de cette
  structure est la norme relativement au morphisme fini et plat $\pi$.
\end{enumerate}
\end{lemme}
\begin{rem}
\label{carre}
Consid{\'e}rons un carr{\'e} 
$
\left(
\begin{array}{c}
\xymatrix{N\ar@{-}[r]\ar@{-}[d]&P\ar@{-}[d]\\
L\ar@{-}[r]&M
}
\end{array}
\right)
$
dans $PIC(X)$ correspondant {\`a} un isomorphisme 
$\phi: L\otimes P \isom M \otimes N$ et soient des isomorphismes 
$\ga : L \isom M$ et $\gb : N \isom P$ tels que le diagramme 
\[
\begin{CD}
 L\otimes P @>{\phi}>> M \otimes N \\
@V{\ga \otimes \text{id}}VV @V{\text{id}\otimes \gb}VV \\
M \otimes P @= M \otimes P
\end{CD}
\]
soit commutatif, alors le diagramme suivant, obtenu par application de
$\det \pi \lst$ l'est aussi:
\[
\begin{CD}
 \det\pi\lst L\otimes \det\pi\lst P @>{\gd_{\phi}}>> 
\det\pi\lst M \otimes \det\pi\lst N \\
@V{\det (\ga )\otimes \text{id}}VV @V{\text{id}\otimes \det (\gb )}VV \\
\det\pi\lst M \otimes \det\pi\lst P @= 
\det\pi\lst M \otimes \det\pi\lst P
\end{CD}
\]
\end{rem}
\subsection{Restriction {\`a} un diviseur effectif}
\label{restriction}
Si $\pi : X\F S$ un morphisme projectif et plat sur $S$ localement
noeth{\'e}rien et $D$ est un diviseur 
relatif sur $X$, si $K$ est un p-cube dans $\cc = PIC(X)$, consid{\'e}rons
le (p+1)-cube 
$A =( \xymatrix{K \ar@{-}[r]_-i &K\otimes \oo (D) })$. On a alors un
isomorphisme canonique 
$r_i: \gd (K\pr ) \isom \gd (K\otimes \oo (D)|_D)$ dans $\dd = PICgr (S)$,
qu'on appellera 
{\em isomorphisme de restriction}. Il est donn{\'e} en appliquant,  pour
chaque sommet 
$L$ de $K$, le foncteur $\gd$ {\`a} la suite exacte:
\[
0 \F L \F L(D) \F L(D)|_D \F 0
\]
\subsubsection{Application au cas des courbes}
\label{cas-des-courbes}
Soit $\pi :X\F S$ un morphisme projectif et plat, {\`a} fibres de dimension
1, sur un sch{\'e}ma localement noeth{\'e}rien $S$. Soient $D$ et $E$ deux
diviseurs de Cartier relatifs effectifs et $L$ un
faisceau inversible sur $X$, consid{\'e}rons alors le carr{\'e} 
$A=K_L(\ox (D),\ox (E))$. Les consid{\'e}rations pr{\'e}c{\'e}dentes
nous donnent un isomorphisme 
$
r_1:\gd(A) \isom \gd(
\xymatrix{
L(D)|_D \ar@{-}[r]_-1&L(D+E)|_D
})
$ et, comme $D$ est fini et plat sur $S$, on obtient un isomorphisme 
\begin{equation}
\label{isom-dim1}
\gd (A)\isom
(\det \pi \lst L(D)|_D)^{-1} \otimes (\det \pi \lst L(D+E)|_D).
\end{equation}
 La section canonique de $\ox(E)$ induit un morphisme 
$\pi \lst L(D)|_D \F \pi \lst L(D+E)|_D$, dont le d{\'e}terminant d{\'e}finit
donc une section $s_L(D,E)$ de $\gd(A)$.
\addtocounter{theo}{1}
\begin{lemme}
\label{indep2}
L'isomorphisme canonique 
$\gd(K_L(\ox (D),\ox (E))) \isom \gd(K_L(\ox (E),\ox (D)))$, donn{\'e} par
les isomorphismes de commutativit{\'e} dans $PICgr(S)$, {\'e}change
les sections $s_L(D,E)$ et $s_L(E,D)$ de ces deux faisceaux.
\end{lemme}
\begin{proof}[Preuve]
D{\'e}crivons l'isomorphisme (\ref{isom-dim1}). Consid{\'e}rons le diagramme
commutatif de suites exactes:
\[
\xymatrix{
&L(E)|_E \ar@{^{(}->}[rr]\ar@{}[rd]|{B}
&&L(D+E)|_{D+E} \ar@{->>}[rr]
&&L(D+E)|_D \\
0\ar[ru]\ar[rr]&&L(D)|_D\ar[ru]\ar@{=}[rr]&&L(D)|_D \ar[ru]\\
&L(E)\ar@{^{(}->}'[r][rr]\ar@{->>}'[u][uu]\ar@{}[rd]|{A}
&&L(D+E)\ar@{->>}'[r][rr]\ar@{->>}'[u][uu]
&&L(D+E)|_D\ar[uu] \\
L\ar@{^{(}->}[rr]\ar[ru]\ar[uu]
&&L(D)\ar@{->>}[rr]\ar[ru]\ar@{->>}[uu]
&&L(D)|_D\ar[uu]\ar[ru] \\
&L\ar@{=}'[r][rr]\ar@{^{(}->}'[u][uu]&&L\ar@{^{(}->}'[u][uu]\\
L\ar@{=}[rr]\ar@{=}[ru]\ar@{=}[uu]&&L\ar@{^{(}->}[uu] \ar@{=}[ru]
}
\]
On en d{\'e}duit que l'isomorphisme (\ref{isom-dim1}) se d{\'e}compose en
$\gd(A)\isom\gd(B)\isom 
(\det \pi \lst L(D)|_D)^{-1} \otimes (\det \pi \lst L(D+E)|_D)$, o{\`u} le
premier isomorphisme est donn{\'e} par les suites exactes verticales et le
second isomorphisme est obtenu en appliquant le foncteur d{\'e}terminant
au diagramme suivant de $\os$-modules localement libres
\[
\begin{CD}
\pi \lst L(E)|_E @>{\gb}>> \pi \lst L(D+E)|_E @. \\
@| @AAA @. \\
\pi \lst L(E)|_E @>>> \pi \lst L(D+E)|_{D+E} @>>> \pi \lst L(D+E)|_D\\
@AAA @AAA @A{\ga}AA\\
0@>>> \pi \lst L(D)|_D @= \pi \lst L(D)|_D
\end{CD}
\]
On traduit l'assertion du lemme en disant que la suite d'isomorphismes
dans $PICgr(S)$
\begin{multline*}
\det( \pi \lst L(D)|_D)\ve \otimes \det( \pi \lst L(D+E)|_D)\isom \\
\det( \pi \lst L(D)|_D)\ve \otimes
\det( \pi \lst L(E)|_E)\ve \otimes\det( \pi \lst L(D+E)|_{D+E}) \isom
\\
\det( \pi \lst L(E)|_E)\ve \otimes\det( \pi \lst L(D)|_D)\ve \otimes
\det( \pi \lst L(D+E)|_{D+E}) \isom \\
\det( \pi \lst L(E)|_E)\ve \otimes \det( \pi \lst L(D+E)|_E)
\end{multline*}
fait correspondre les sections $\det(\ga)$ et  $\det(\gb)$ des deux
termes extr{\`e}mes, ce qui provient des propri{\'e}t{\'e}s du
d{\'e}terminant des $\os$-modules localement libres.
\end{proof}
\addtocounter{ptheo}{2}
\begin{ptheo}
\addtocounter{subsection}{1}
\label{th-prin}
Soit $S$ un sch{\'e}ma localement noeth{\'e}rien. 
Pour tout morphisme projectif et plat $\pi : X \F S$ {\`a} fibres de
dimension n, il existe une structure du (n+2)-cube 
canonique sur le foncteur 
$\gdp : PIC(X) \F PIC (S)$ telle que:
\begin{enumerate}
\item Si $\pi$ est un morphisme fini, la structure du carr{\'e}
  correspondante est simplement celle donn{\'e}e par la norme
  (\ref{norme}).
\item Si $D$ est un diviseur relatif sur $X$, les structures du
  (n+2)-cube sur $\gdp _X: PIC(X)\F PIC(S)$ et du (n+1)-cube sur 
$\gdp _D: PIC(D) \F PIC(S)$ sont
  compatibles aux isomorphismes de restriction.
\end{enumerate}
\end{ptheo}
\begin{proof}[Preuve]
La preuve s'effectue par r{\'e}currence sur $n$.

Le cas de la dimension 0 a {\'e}t{\'e} trait{\'e} dans (\ref{norme}). 

Soit n un entier strictement positif, supposons que, pour tout
 morphisme  projectif et plat $Y\F S$ {\`a} fibres de dimension $p<n$,
 on sache 
 construire une structure du (p+2)-cube sur $\gdp_Y$ verifiant les
 propri{\'e}t{\'e}s 1 et 2 du th{\'e}or{\`e}me. 
Soit $\pi :X\F S$ projectif et plat {\`a} fibres de
 dimension n, construisons une structure du (n+2)-cube sur
 $\gdp_X$. D'apr{\`e}s la proposition (\ref{cubebis}), on doit donc
 construire, pour tout $(n+2)$-cube $A$ dans $PIC(X)$, des
 isomorphismes 
$
s_{A,i} :
 \gd ((\phi \pr_i)\ust A) 
\isom  
\gd ((\phi \sec_i)\ust A
$ 
v{\'e}rifiant des propri{\'e}t{\'e}s de sym{\'e}trie et de
 compatibilit{\'e} aux recollement de $(n+2)$-cubes. Cette
 construction, qui  occupe les paragraphes suivants sera
 d{\'e}coup{\'e}e de la fa\c{c}on suivante: 
 \begin{itemize}
 \item Si un $(n+2)$-cube $A$ poss{\`e}de une ar{\^e}te $L_j$ ayant
   une section 
   r{\'e}guli{\`e}re $\gs_j$, on construit un isomorphisme 
$
s_{A,i,\gs_j} :
 \gd ((\phi \pr_i)\ust A) 
\isom  
\gd ((\phi \sec_i)\ust A
$
   d{\'e}pendant du choix de $\gs_j$.
\item On s'affranchit de la d{\'e}pendance en la section $\gs_j$, sous
  l'hypoth{\`e}se que $A$ poss{\`e}de deux ar{\^e}tes $L_j$ et $L_k$
  suffisamment 
  positives (on dira dans ce cas que $A$ est suffisamment positif dans
les directions $i$ et $j$).

\item On {\'e}limine enfin cette hypoth{\`e}se de positivit{\'e};
 \end{itemize}
\subsection{Construction et propri{\'e}t{\'e}s de $s_{A,i,\gs_j}$ }
\label{constr-avec-div}
\begin{constr} \label{constr1}
Une section r{\'e}guli{\`e}re 
$\gs_j$ de la $j$-i{\`e}me ar{\^e}te $L_j$ de $A$ d{\'e}finit un diviseur de
Cartier relatif $D$ et $A$ est 
isomorphe {\`a} un $(n+2)$-cube
 $( \xymatrix{K \ar@{-}[r]_-j& K\otimes \ox (D)})$. On dispose donc,
 d'apr{\`e}s 
 (\ref{restriction}),  
 d'isomorphismes de restriction
$
 \gd ((\phi \pr  _i)\ust A) \isom \gd ((\phi \pr _i)\ust K(D)|_D)
$
et 
$
 \gd ((\phi \sec  _i)\ust A) \isom \gd ((\phi \sec _i)\ust K(D)|_D)
$,
provenant du diagramme de suites exactes
\[
\begin{CD}
0 @>>> (\phi \pr _i)\ust K  @>>>
 (\phi \pr _i)\ust K(D)  @>>>
  (\phi \pr _i)\ust K(D)|_D  @>>> 0
\end{CD}
\]
et du  diagramme analogue faisant intervenir $(\phi \sec  _i)\ust K$. 
L'hypoth{\`e}se de r{\'e}currence, appliqu{\'e}e au $(n+1)$-cube
$K(D)|_D$ sur le 
sch{\'e}ma relatif $D/S$ de dimension $n-1$ entra{\^\i}ne l'existence d'un
isomorphisme  
\[
s_{K(D)|_D,i} :
 \gd ((\phi \pr _i)\ust K(D)|_D ) 
\isom  
\gd ((\phi \sec _i)\ust K(D)|_D).
\]
En composant $s_{K(D)|_D,i}$ avec les isomorphismes pr{\'e}c{\'e}dents, on
obtient $s_{A,i,\gs _j}$.
\end{constr}
\addtocounter{subsubsection}{1}
\subsubsection{Comportement de $s_{A,i,\gs _j}$ par recollement.}
\label{recol_div}
Si deux $(n+2)$-cubes $A$ et $B$ sont  recollables le long de leurs
i-{\`e}me face, ils ont alors m{\^e}me j-i{\`e}me ar{\^e}te $L_j$. Si
$L_i$ poss{\`e}de 
une section r{\'e}guli{\`e}re $\gs_j$ d{\'e}finissant un diviseur
relatif $D$, on 
peut alors {\'e}crire $A\simeq(\xymatrix{K\ar@{-}[r]_-j &K\otimes\oo(D)})$ 
et $B=\simeq(\xymatrix{K\pr \ar@{-}[r]_-j &K\pr\otimes\oo(D)})$, avec 
 $K=(  \xymatrix{K_1 \ar@{-}[r]_-i &K_2 })$ et
 $K\pr =(  \xymatrix{K_2 \ar@{-}[r]_-i &K_3 })$. En {\'e}crivant le diagramme
 suivant, dont chaque ligne est exacte:
\[
\xymatrix{
0 \ar[r] & K_3 \ar@{-}[d]  \ar[r]
\ar@{}[dr]|{B} &
  K_3(D) \ar@{-}[d]  \ar[r]
 &  K_3(D)|_D  \ar@{-}[d]  \ar[r]
&0 \\
0 \ar[r] & K_2  \ar@{-}[d]_-i  \ar[r]
\ar@{}[dr]|{A}&
  K_2(D)  \ar@{-}[d] \ar[r]
 & K_2(D)|_D  \ar@{-}[d]  \ar[r]
&0 \\
0 \ar[r] & K_1    \ar[r]_-j
&  K_1(D)   \ar[r]
 &  K_1(D)|_D   \ar[r]
&0
}
\]
et en utilisant que, d'apr{\`e}s l'hypoth{\`e}se de r{\'e}currence, 
\[
s_{(K(D)|_D \ast _i K\pr(D)|_D),i}
 = 
s_{ K\pr(D)|_D,i} \circ s_{K(D)|_D ,i}\; ,
\]
on obtient:
\[
s_{A \ast _i B,i,\gs _j}
 = 
s_{B,i,D} \circ s_{A,i,\gs _j}.
\]
\begin{lemme}[Lien entre $s_{A,i,\gs _j}$ et  $s_{A,i,\gs _k}$]
\label{indep} 
Supposons que les ar{\^e}tes $L_j$ et $L_k$ du $n+2$-cube $A$ poss{\`e}dent
des sections r{\'e}guli{\`e}res $\gs _j$ et $\gs _k$ et soit $i\neq j,k$. Si 
l'une des deux hypoth{\`e}ses suivantes est v{\'e}rifi{\'e}e:
\begin{enumerate}
\item $X/S$ est {\`a} fibres de dimension 1 et $Z(\gs_j)\cap
  Z(\gs_k)=\emptyset$.
\item $X/S$ est {\`a} fibres de dimension $n>1$ et les suites 
 $(\gs _j,\gs _k)$ et  $(\gs _k,\gs _j)$ sont $\pi$-r{\'e}guli{\`e}res
\end{enumerate}
On a alors $s_{A,i,\gs _j}= s_{A,i,\gs _k}$.
\end{lemme}
\begin{proof}[Preuve]
Soient $D$ et $E$ les diviseurs de Cartier relatifs effectifs d{\'e}finis par 
 $s_{A,i,\gs _j}$ et  $s_{A,i,\gs _k}$, on {\'e}crit alors  $A$ sous la forme 
$
 A=\left(
 \begin{array}{c}
 \xymatrix{
 K(E) \ar@{-}[r]\ar@{-}[d]_-k &K (D+E)\ar@{-}[d] \\
 K \ar@{-}[r]_-j &K (D)
 }
\end{array}
\right)
$, o{\`u} $K$ est un $n$-cube. Ecrivons le diagramme suivant,
dont les lignes et les colonnes sont des suites exactes
courtes, et qui 
relie le $(n+2)$-cube $A$ sur X aux $(n+1)$-cubes $A\pr$ et $A\sec$
sur $E$ et $D$ et au $n$-cube $A ^{\prime\prime\prime}$ sur $D\cap E$:
\begin{equation}
\label{double-restriction}
\xymatrix{A\sec &
K(E)|_E \ar@{^{(}->}[r]_-j&K(D+E)|_E\ar@{->>}[r]
&K(E)|_{D\cap E}
&A ^{\prime\prime\prime} \\
A&
K(E)\ar@{^{(}->}[r]\ar@{->>}[u] &K(D+E)\ar@{->>}[r]\ar@{->>}[u]
&K(D+E)|_D\ar@{->>}[u]
&A\pr\\
&K\ar@{^{(}->}[r]_-j\ar@{^{(}->}[u]_-k &
\hspace{4mm} K(D) \hspace{4mm}
\ar@{^{(}->}[u]\ar@{->>}[r]&K(D)|_D\ar@{^{(}->}[u]_-k
\save "2,2"."3,3"*[F-]\frm{}\ar@{.}"2,1"\restore
\save "1,2"."1,3"*[F-]\frm{}\ar@{.}"1,1"\restore
\save "2,4"."3,4"*[F-]\frm{}\ar@{.}"2,5"\restore 
\save "1,4".*[F-]\frm{}\ar@{.}"1,5"\restore 
}
\end{equation}
{\em Cas 1}:\\
$K$ est ici un 1-cube $(\xymatrix{L\ar@{-}[r]&M})$ et le diagramme
(\ref{double-restriction}) se r{\'e}duit alors {\`a}:
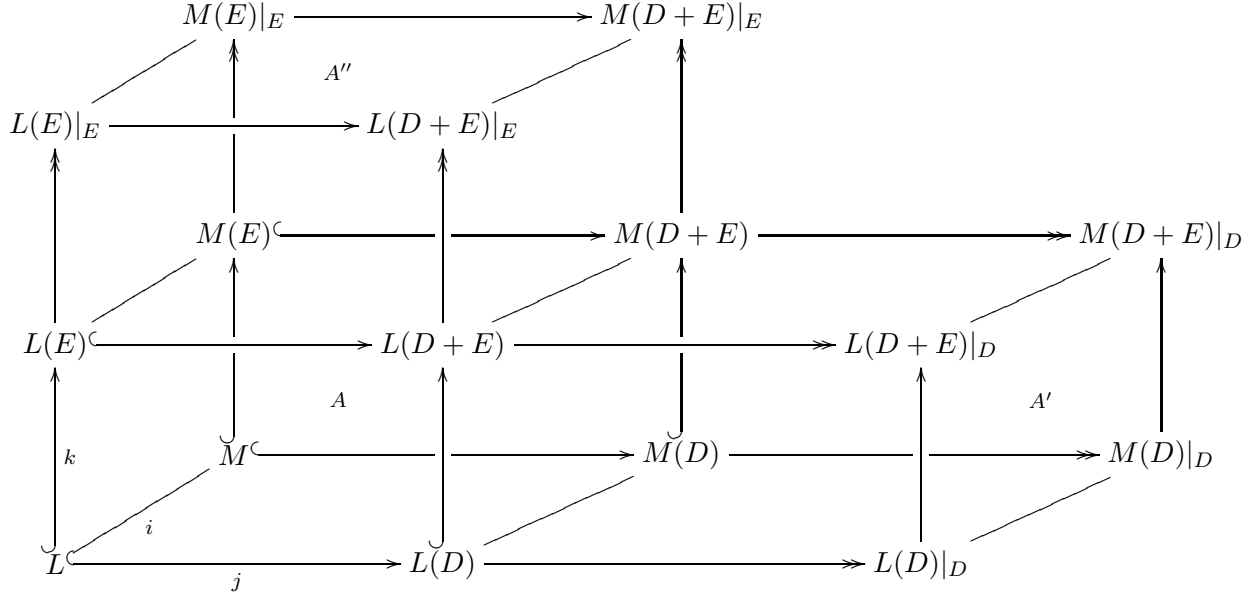
\[
\xymatrix{
&M(E)|_E\ar[rr] \ar @{} [dr] |{A\sec}
 && M(D+E)|_E \\
L(E)|_E\ar[rr]\ar@{-}[ur]
 && L(D+E)|_E\ar@{-}[ur]&\\
&M(E)\ar@{^{(}->}'[r][rr]\ar@{->>}'[u][uu]
 && M(D+E)\ar@{->>}[uu]\ar@{->>}[rr]&&M(D+E)|_D \\
L(E)\ar@{^{(}->}[rr]\ar@{->>}[uu]\ar@{-}[ur]
 && L(D+E)\ar@{->>}[uu]\ar@{-}[ur]\ar@{->>}[rr]
&&L(D+E)|_D\ar@{-}[ur]\ar @{} [rd] |{A\pr}&\\
&M\ar@{^{(}->}'[r][rr]\ar@{^{(}->}'[u][uu] \ar @{} [ur] |{A}
 && M(D)\ar@{->>}'[r][rr]\ar@{^{(}->}'[u][uu]
&&M(D)|_D\ar[uu]  \\
L\ar@{^{(}->}[rr]_-j\ar@{-}[ur]_-i\ar@{^{(}->}[uu]_-k &&
 L(D)\ar@{-}[ur]\ar@{->>}[rr]\ar@{^{(}->}[uu]&&L(D)|_D\ar[uu]\ar@{-}[ur]&
}
\]
o{\`u} $\{i,j,k\} =\{1,2,3\}$. On peut alors {\'e}crire un diagramme
d'isomorphismes: 
\[
\xymatrix{
\gd(\fipr_i(A\sec)) \ar[dd]^-{s_{A\sec,i}}
&\gd(\fipr_i(A)) \ar[l] \ar[r]
&\gd(\fipr_i(A\pr))  \ar[dd]_-{s_{A\pr,i}}\\
& {\os}\ar[ru] \ar[rd]\ar[lu] \ar[ld] \\
\gd(\fisec_i(A\sec)) 
&\gd(\fisec_i(A)) \ar[l] \ar[r]
&\gd(\fisec_i(A\pr))
}
\]
dont les lignes sup{\'e}rieures et inf{\'e}rieures proviennent respectivement
des faces avant et arri{\`e}res du diagramme pr{\'e}c{\'e}dent et les
fl{\`e}ches 
obliques sont les isomorphismes d{\'e}crits  dans
(\ref{cas-des-courbes}). Il s'agit de montrer 
que le trac{\'e} ext{\'e}rieur est commutatif, ce qui provient du fait que les
deux triangles lat{\'e}raux sont commutatifs par la remarque (\ref{carre})
et que les deux
autres triangles sont commutatis par le lemme (\ref{indep}).\\
{\em Cas 2}\\
Le diagramme (\ref{double-restriction}) induit un diagramme commutatif
d'isomorphismes 
\[
\begin{CD}
\gd (A\sec )  @>{\gd}>>\gd (A ^{\prime\prime\prime} )\\
@A{\gc}AA @A{\gb}AA \\
\gd (A) @>{\ga}>>\gd (A\pr )
\end{CD}
\]
Par construction, $\ga $ identifie $s_{A\pr ,i}$ {\`a} $s_{A,i,\gs_j}$ et 
 $\gc $ identifie $s_{A\sec ,i}$ {\`a} $s_{A,i,\gs_k}$. Par l'hypoth{\`e}se de
 r{\'e}currence appliqu{\'e}e {\`a} $\gd _D$, $\gd _E$ et  $\gd
 _{D\cap E}$, les 
 isomorphismes $\gb$ et  $\gd$ identifient respectivement $s_{A\pr ,i}$
 et $s_{A\sec ,i}$ {\`a}  $s_{A ^{\prime \prime\prime} ,i}$, ce qui prouve
 l'assertion dans le cas 2.
\end{proof}
\subsubsection{Lien entre $s_{A,i,\gs_j}$ et $s_{A,k,\gs_j}$. }
\label{signe_div}
Soient $A$ un $(n+2)$-cube sur $X$ dont la $j$-i{\`e}me ar{\^e}te poss{\`e}de 
 une section $\pi$-r{\'e}guli{\`e}re $\gs_j$ et soent $i,k$ deux indices
 distincts et distincts de $j$, montrons que 
les trivialisations de $\gd (A)$ induites par  $s_{A,i,\gs_j}$ et
$s_{A,k,\gs_j}$ diff{\`e}rent d'un signe {\'e}gal {\`a} $\ge _{ik}
 (A)$ (introduit 
dans la proposition (\ref{cubebis})).\\
 A cet effet, {\'e}crivons $K$ sous la forme 
$\left(
\begin{array}{c}
\xymatrix{
 K_2 \ar@{-}[r]\ar@{-}[d]_-k &K_3\ar@{-}[d] \\
 K_0 \ar@{-}[r]_-i &K_1
 }
\end{array}
\right)
$
et notons $D$ le diviseur de Cartier relatif d{\'e}fini par $\gs_j$.
Consid{\'e}rons le diagramme suivant dont les lignes sont des suites
exactes courtes:
\[
\xymatrix{
&K_3\ar@{^{(}->}[rr]\ar@{-}'[d][dd]\ar@{-}[dl]&&
K_3(D)\ar@{->>}[rr]\ar@{-}'[d][dd]\ar@{-}[dl]
&&K_3(D)|_D \ar@{-}[dl]\ar@{-}[dd]\\
K_2\ar@{^{(}->}[rr] \ar@{-}[dd]&&
K_2(D)\ar@{->>}[rr]\ar@{-}[dd]&&K_2(D)|_D\ar@{-}[dd]& \\
&K_1\ar@{^{(}->}'[r][rr]\ar@{-}[dl]&&K_1(D)\ar@{-}[dl]\ar@{->>}'[r][rr]
&&K_1(D)|_D\ar@{-}[dl] \\
K_0\ar@{^{(}->}[rr]&&K_0(D)\ar@{->>}[rr]&&K_0(D)|_D&
}
\] 
le cube de gauche de ce diagramme est pr{\'e}cis{\'e}ment $A$ et les suites
exactes donnent un isomorphisme entre $\gd (A) $ et 
$\gd \left(
\begin{array}{c}
 \xymatrix{
 K_2(D)|_D \ar@{-}[r]\ar@{-}[d]_-j &K_3(D)|_D\ar@{-}[d] \\
 K_0(D)|_D \ar@{-}[r]_-i &K_1(D)|_D
 }
\end{array}
\right)
$.
Cet isomorphisme identifie 
les trivialisations de $\gd (A)$ induites par $s_{A,i,\gs_j}$ et
$s_{A,k,\gs_j}$  aux trivialisations de $\gd (K(D)|_D)$
induites par $s_{K(D)|_D,i}$ et $s_{K(D)|_D,k}$. Or, par l'hypoth{\`e}se de
r{\'e}currence, interpr{\'e}t{\'e}e {\`a} la lumi{\`e}re de la proposition
(\ref{cubebis}), ces deux 
trivialisations diff{\`e}rent d'un signe $\ge _{ik}(K(D)|_D)$. On conclut
en affirmant que $\ge _{ik}(A) = \ge _{ik}(K(D)|_D)$. En effet, par
l'additivit{\'e} de la caract{\'e}ristique d'Euler-Poincar{\'e}, on a:
\[
\chi_{D/S} (K_0 (D)|_D) =\chi_{X/S} (K_0(D))-\chi_{X/S} (K_0)
=\chi_{X/S}(\xymatrix{K_0 \ar@{-}[r] &K_0 (D)})
\]
et donc:
\[
\ge _{ik}(A) = (-1)^{\chi_{X/S}(K_0)-\chi_{X/S}(K_0 (D))} 
=(-1)^{\chi_{D/S}(K_0(D)|_D)} 
= \ge_{ik} (K(D)|_D)
\] 
\subsection{Elimination des hypoth{\`e}ses de diviseurs effectifs}
Pour tout $(n+2)$-cube $A$  suffisamment positif dans les directions $i$
et $j$ et tout indice $i\neq j,k$, on construit dans
cette section un isomorphisme 
\[
s_{A,i}: \gd (\fipr_i A) \isom  \gd (\fisec_i A),
\]
fonctoriel en les isomorphismes de $(n+2)$-cubes et compatible aux
empilements de cubes dans la direction $i$ (conditions 1 et 2 de la
proposition (\ref{cubebis})).
\subsubsection{Premi{\`e}re r{\'e}duction}
\label{reduction}
Soit $A$ un $(n+2)$-cube suffisamment positif dans les directions $i$
et $j$. $A$ est alors isomorphe un $(n+2)$-cube 
$
\left(
\begin{array}{c}
 \xymatrix{
 K\otimes M\ar@{-}[r]\ar@{-}[d]_-k &K\otimes (L\otimes M)\ar@{-}[d]\\
 K         \ar@{-}[r]_-j &K\otimes L
 }
\end{array}
\right)
$ 
avec $L,M \gg 0$.
Consid{\'e}rons les fibr{\'e}s projectifs $P_L$ et $P_M$ sur $S$ et effectuons
le changement de base  
\[
\begin{CD}
X_{P_L \times P_M} @>g>> X\\
@V{\pi}VV @V{\pi}VV \\
P_L \times P_M @>f>> S
\end{CD}
\]
Notons encore $K$, $L$ et $M$ les images r{\'e}ciproques de $K$, $L$ et
$M$ par $g$ 
et introduisons les faisceaux inversibles 
$L\pr = L\otimes \pi_{P_L}\ust \oo_{P_L}(1)$ 
et 
$M\pr = M\otimes \pi_{P_M}\ust \oo_{P_M}(1)$ et consid{\'e}rons le
$(n+2)$-cube 
$A\pr =
\left(
\begin{array}{c}
 \xymatrix{
 K\otimes M\pr\ar@{-}[r]\ar@{-}[d]_-k &K\otimes (L\pr\otimes M\pr)\ar@{-}[d]\\
 K         \ar@{-}[r]_-j &K\otimes L\pr
 }
\end{array}
\right)
$ 
sur $X_{P_L \times P_M}$. On a des isomorphismes canoniques:
$
\gd(K\otimes L\pr) \isom f\ust\gd(K\otimes L) 
\otimes 
(\oo_{P_L}(1))^{\chi_{X/S}(K\otimes L)}
$, 
$
\gd(K\otimes M\pr) \isom f\ust\gd(K\otimes M) 
\otimes 
(\oo_{P_M}(1))^{\chi_{X/S}(K\otimes M)}
$ et  
$
\gd(K\otimes L\pr\otimes M\pr) 
\isom f\ust\gd(K\otimes L\otimes M) 
\otimes 
(\oo_{P_L}(1)\otimes\oo_{P_M}(1))^{\chi_{X/S}(K\otimes L\otimes M)}
$. L'{\'e}galit{\'e} 
$\chi_{X/S}(K\otimes L\otimes M)=
\chi_{X/S}(K\otimes L)=\chi_{X/S}(K\otimes M)$, d{\'e}duite de
(\ref{caracteristique2}),
entraine alors l'existence d'un isomorphisme canonique 
\[
\gd(A\pr) \isom f\ust \gd(A)
\]
 Comme on a 
$f\lst f\ust \oo _{P_L \times _S P_M} = \oo _S$, on en d{\'e}duit que la
donn{\'e}e d'un isomorphisme 
$s_{A,i}: \gd (\fipr_i A) \isom  \gd (\fisec_i A)$ sur $S$ est
{\'e}quivalente {\`a} la donn{\'e}e d'un  isomorphisme 
$s_{A\pr,i}: \gd (\fipr_i A\pr) \isom  \gd (\fisec_i A\pr)$ sur 
$P_L\times _S P_M$.
\subsubsection{Construction de $s_{A\pr,i}$ sur des ouverts}
Consid{\'e}rons les ouverts 
 $U_L = P_L \setminus Z_L$ et  $U_M = P_M \setminus Z_M$
introduits en (\ref{div2}). 
Sur $X_{U_L \times P_M}$, on dispose d'un diviseur relatif $D_L$ et
d'un isomorphisme canonique 
$L\pr \isom \oo (D_L)$. Le cube $A\pr$ est donc
canoniquement isomorphe sur  $X_{U_L \times P_M}$ {\`a} un cube de la forme
$(\xymatrix{ K_0 \ar@{-}[r] & K_0 (D_L)})$. En appliquant la
construction (\ref{constr1}) {\`a} cette situation, on obtient un
isomorphisme canonique $s\pr :\gd (\fipr_i A) 
\isom \gd (\fisec_i A)$ d{\'e}fini sur $U_L\times _S P_M$. La m{\^e}me
construction, effectu{\'e}e en intervertissant les r{\^o}les de $L$ et $M$,
nous donne un autre isomorphisme  $s\sec$ entre les m{\^e}mes faisceaux et
d{\'e}fini sur  $P_L\times _S U_M$.
\addtocounter{theo}{2}
\addtocounter{subsubsection}{1}
\begin{lemme}
\label{coincid}
$s\pr$ et $s\sec$ coincident sur $U_L \times _S U_M$.
\end{lemme}
\addtocounter{subsubsection}{1}
\begin{proof}[Preuve]
Soit $p$ la projection $U_L \times _S U_M \F U_M$. Le lemme
(\ref{genericite2}) montre l'existence d'un ferm{\'e} $Z$ de  $U_L \times
_S U_M$ tel que
\begin{enumerate}
\item $\forall u \in U_M, \forall x \in Z_u, 
\text{Prof}(\oo_{p ^{-1}(u),x}) \geq 1$.
\item Au dessus de $V=(U_L \times _S U_M) \setminus Z$, les diviseurs
  $D_L$ et $D_M$ sont en position d'intersection compl{\`e}te et 
$(D_L \cap D_M)_V$ est plat sur $V$.
\end{enumerate}
Au dessus de $V$, $s\pr$ et $s\sec$ coincident d'apr{\`e}s
(\ref{indep}). On en d{\'e}duit donc, d'apr{\`e}s (\cite{EGA4},19.9.8), en
utilisant la 
propri{\'e}t{\'e} 1 que $s\pr$ et $s\sec$ coincident sur $U_L \times _S U_M$.
\end{proof}
\subsubsection{Extension de $s\pr$ et $s\sec$ {\`a} $P_L \times _S P_M$.}
Consid{\'e}rons le ferm{\'e} $Y=Z_L \times _S Z_M$ du $s$-sch{\'e}ma
$P_L\times _S P_M$. Le lemme (\ref{genericite1}) entra{\^\i}ne que
pour tout $s\in S$ et 
tout $y\in Y_s$, on a $\text{Prof}(\oo_{(P_L \times _S P_M)_s,y}) \geq
2$. les isomorphismes $s\pr$ et $s\sec$ sont d{\'e}finis respectivement
sur $U_L \times _S P_M$ et $P_L \times _S U_M$ et coincident sur 
$U_L \times _S U_M$. Ils d{\'e}finissent donc un isomorphisme 
$s: \gd ((\phi_i\pr)\ust A) \isom  \gd ((\phi_i\sec )   \ust A)$ sur
l'ouvert $U= (P_L \times _S P_M)\setminus Y$. Par
(\cite{EGA4},19.9.8), $s$ se prolonge alors (de mani{\`e}re unique) {\`a} 
$P_L \times _S P_M$ en un isomorphisme 
$s_{A\pr,i}^{j,k} : \gd (\fipr_i A\pr) \isom \gd (\fisec_i
A\pr)$
 d{\'e}pendant  a priori du choix
des directions $j$ et $k$.\\
Par la r{\'e}duction (\ref{reduction}), on d{\'e}duit de cette construction un
isomorphisme 
\[
s_{A\pr,i}^{j,k} : \gd (\fipr_i A\pr) \isom \gd (\fisec_i
A\pr)
\]
d{\'e}fini maintenant sur $S$. Il reste {\`a} montrer le
\addtocounter{theo}{1}
\begin{lemme}
$s_{A,i}^{j,k}$ est ind{\'e}pendant du choix des directions $j$ et $k$.
\end{lemme}
\begin{proof}[Preuve]
Supposons en effet que $A$ soit suffisamment positif dans trois
directions diff{\'e}rentes $j,k,l \neq i$. $A$ s'{\'e}crit donc de trois
fa\c{c}ons diff{\'e}rentes 
$(\xymatrix{K_j \ar@{-}[r]_-j & K_j \otimes L})$, 
$(\xymatrix{K_k \ar@{-}[r]_-k & K_k \otimes M})$ et 
$(\xymatrix{K_l \ar@{-}[r]_-l & K_l \otimes N})$, 
avec $L,M,N \gg 0$. Montrons alors $s_{A,i}^{j,k}=s_{A,i}^{k,l}$.\\
Il suffit de le montrer apr{\`e}s changement de base par $P_M \F S$, et
donc apr{\`e}s changement de base par $U_M \F S$ (par le m{\^e}me argument que
dans la preuve du lemme (\ref{coincid}) ). Mais, par construction, au
dessus de $U_M$ les 
isomorphismes  $s_{A,i}^{j,k}$ et $s_{A,i}^{k,l}$ proviennent tous
deux de la structure du cube associ{\'e}e au sch{\'e}ma relatif $D_M/S$, et
sont donc {\'e}gaux.
\end{proof}
\addtocounter{subsubsection}{1}
\subsubsection{Propri{\'e}t{\'e}s des isomorphismes ainsi construits} On a
donc construit, pour tout $(n+2)$-cube $A$, qui 
est suffisamment positif dans deux directions et pour toute direction
$i$ diff{\'e}rente des pr{\'e}c{\'e}dentes, un isomorphisme
\[
s_{A,i} : \gd (\fipr_i A) \isom \gd (\fisec_i
A)
\]
les $s_{A,i}$ v{\'e}rifient les propri{\'e}t{\'e}s suivantes:
\begin{enumerate}
\item Si les deux  $(n+2)$-cubes $A$ et $B$ sont suffisamment positifs
  dans deux directions diff{\'e}rentes de $i$ et ont leur i-{\`e}me face en
  commun, alors:
\[
s_{B,i} \circ s_{A,i} = s_{A\ast _i B,i}
\]
\item Si $s_{A,i}$ et $s_{A,j}$ sont tous les deux d{\'e}finis, alors les
  trivialisations de $\gd (A)$ qu'elles induisent diff{\`e}rent d'un signe
  {\'e}gal {\`a} $\ge _{ij} (A)$.
\item Si $A$ est suffisamment positif dans deux directions diff{\'e}rentes
  de $i$, pour toute permutation $\gs \in S_n$ on a:
\[
s_{\gs \ust A,i} = s_{A, \gs (i)}
\]
\end{enumerate}
Pour montrer les  propri{\'e}t{\'e}s 1 et 2, on choisit une direction $k$
dans laquelle $A$ est suffisamment positif, et on {\'e}crit
 $A=(\xymatrix{K \ar@{-}[r]_-k &K \otimes L})$ avec $L\gg 0$ et on
 effectue,  comme pr{\'e}c{\'e}demment, 
un changement de base $U_L \subset P_L \F S$, ce qui permet de disposer
d'un diviseur effectif $D_L$. La propri{\'e}t{\'e} 1 provient alors de
la propri{\'e}t{\'e} de 
recollement (\ref{recol_div}) et la propri{\'e}t{\'e} 2 provient du lien
entre $s_{A,i,D}$ 
et  $s_{A,k,D}$ (\ref{signe_div}), montr{\'e}s tous les deux sous
l'hypoth{\`e}se de l'existence d'un diviseur effectif.\\
Pour la propri{\'e}t{\'e} 3, il
suffit de montrer l'{\'e}galit{\'e} $s_{\gs \ust A,i} = s_{A, \gs (i)}$ dans
le cas o{\`u}  $\gs$ est une
transposition $\gs _{rs}$ et o{\`u} $A$ est suffisamment positif dans au
moins une direction $k$ diff{\'e}rente de $r$ et de $s$. On {\'e}crit alors 
 $A=(\xymatrix{K \ar@{-}[r]_-k &K \otimes L})$ avec $L\gg 0$ et on
 effectue le  changement de base $U_L \subset P_L \F S$. L'{\'e}galit{\'e}
 provient alors de l'{\'e}galit{\'e} analogue pour le $(n+1)$-cube $K\otimes L
 |_{D_L}$, qui est v{\'e}rifi{\'e}e d'apr{\`e}s l'hypoth{\`e}se de
 r{\'e}currence. 
\subsection{Elimination des hypoth{\`e}ses de positivit{\'e}.}
On veut {\'e}tendre ici la d{\'e}finition de $s_{A,i}$ {\`a} un $(n+2)$-cube
quelconque $A$ dans $PIC(X)$.
\begin{rem}
\label{compar-sign}
Si $A$ et $B$ sont deux $(n+2)$-cubes dans $PIC(X)$, recollables le
long de leur j-i{\`e}me face (notons $C= A\ast _j B$) et tels que
$s_{A,i}$ et  $s_{B,i}$ sont bien 
d{\'e}finis. On peut alors d{\'e}finir un isomorphisme
$s_{A,i}\otimes s_{B,i}:\gd (\fipr_i C) \isom 
\gd (\fisec_i C)$ par:
\[
\gd (\fipr_i C) 
\isom 
\gd (\fipr_i A) \otimes \gd (\fipr_i B)
\F
\gd (\fisec_i A) \otimes \gd (\fisec_i B)
\isom 
\gd (\fisec_i C) 
\]
Notons que si $s_{C,i}$ est aussi d{\'e}fini, on a la relation:
\[
s_{C,i} =\ge _{ij}(A) s_{A,i}\otimes s_{B,i}\; ,
\]
r{\'e}sultant de la comparaison entre $s_{A,i}$ et $s_{A,j}$ et de la
relation $s_{C,j} =s_{A,j}\otimes s_{B,j}$. 
\end{rem}
On peut alors {\'e}noncer le
\begin{lemme}
Soit $i$ un indice compris entre 1 et $n+2$ et soit  $A$ un
$(n+2)$-cube suffisamment positif dans {\em une} direction $j\neq i$.
\begin{enumerate}
\item Pour toute direction $k\neq i,j$, il existe un $(n+2)$-cube $B$
  recollable avec $A$ dans la direction $k$ tel que $B$ 
et $A\ast _k B$ soient tous les deux suffisamment positifs dans la
direction $k$.
\item Soient $k$ et $l$ deux directions eventuellement {\'e}gales mais
  distinctes de $i$ et $j$ et soient $B$ (resp. $B\pr$) deux $(n+2)$-cubes
  recollables avec $A$ dans les directions $k$ (resp. $l$) 
  tels que $B$ et $A\ast _k B$ sont suffisament positifs dans la
  direction $k$ et $B\pr$ et $A\ast _k B\pr$ sont suffisamment positifs dans la
  direction $l$. Alors les deux isomorphismes:  
\[
  s_{A\ast_k B ,i} \otimes  s_{B ,i}^{-1}
\; \text{et} \; 
  s_{A\ast_l B\pr ,i}\otimes  s_{B\pr ,i}^{-1}: 
\gd (\fipr_i A) \isom \gd (\fisec_i A) 
\]
sont {\'e}gaux.
\end{enumerate}
\end{lemme}
\begin{proof}[Preuve]
Le premier point r{\'e}sulte du fait que pour tout faisceau inversible $L$
sur $X$, on peut trouver un faisceau inversible $M\gg 0$ tel que 
$L\otimes M \gg 0$.\\
Montrons d'abord le second point quand $k=l$. Dans ce cas, on peut
trouver des $(n+2)$-cubes $C$ et $C\pr$, suffisamment positifs dans la
direction $k$ , tels que 
$A\ast _k B \ast _k C = A\ast _k B\pr \ast _k C\pr$ et $ B \ast _k C$,
$B\pr \ast _k C\pr$ et $A\ast _k B \ast _k C$ sont suffisamment
positifs dans la direction $k$. On a alors en utilisant la remarque
(\ref{compar-sign}: 
\[
s_{A\ast _kB,i} \otimes s_{B,i}^{-1} =
s_{A\ast _kB \ast _k C,i}
\otimes s_{B \ast _k C,i}^{-1} =
s_{A\ast _kB\pr \ast _k C\pr ,i}
\otimes s_{B\pr \ast _k C\pr ,i}^{-1} =
s_{A\ast _kB\pr ,i} \otimes s_{B\pr ,i}^{-1}
\]
Dans le cas  o{\`u} $k\neq l$, on introduit un $(n+2)$-cube $C$
d{\'e}termin{\'e} 
uniquement par les conditions suivantes:
$C$ est recollable avec $B$ dans la direction $l$ et avec $B\pr$ dans la
direction $k$, comme le d{\'e}crit le diagramme suivant:
\begin{center}
\begin{tabular}{cc}
$l \uparrow$ &\begin{tabular}{|c|c|}
\hline
$B\pr$&$C$\\ \hline
$A$&$B$ \\ \hline
\end{tabular}
\\
&$\xrightarrow{k}$
\end{tabular}
\end{center}
En notant $D$ le cube total, on {\'e}crit alors:
\[
s_{A\ast _kB,i} \otimes s_{B,i}^{-1} =
s_{C,i} \otimes s_{B\pr \ast _kC,i}^{-1} \otimes s_{B,i}^{-1}=
s_{C,i} \otimes s_{B \ast _lC,i}^{-1} \otimes s_{B\pr ,i}^{-1}=
s_{A\ast _lB\pr ,i} \otimes s_{B\pr ,i}^{-1}
\]
\end{proof}
Ce lemme permet donc de d{\'e}finir 
$s_{A,i}:\gd (\fipr_i A) \isom \gd (\fisec_i A)$
d{\'e}s que $A$ est suffisamment positif dans {\em une} direction $j\neq
i$ par la formule:
\[
s_{A,i}=s_{B ,i}^{-1} \circ s_{A\ast_k B ,i}
\]
En r{\'e}it{\'e}rant cette argument, on construit $s_{A,i}$ pour tout
$(n+2)$-cube $A$ dans $PIC(X)$, sans hypoth{\`e}se de positivit{\'e} sur $A$.
\begin{lemme}
Les isomorphismes $s_{A,i}$ ainsi construits v{\'e}rifient les
propri{\'e}t{\'e}s de la proposition (\ref{cubebis}).
\end{lemme}
Ce qui conclut la preuve du th{\'e}or{\`e}me principal.
\end{proof}
\section{Fibr{\'e} d'intersection et r{\'e}sultant}
\subsection{Fibr{\'e} d'intersection}
Pour tout morphisme projectif et plat $\pi:X \F S$, {\`a} fibres de
dimension $n$, sur un sch{\'e}ma $S$ localement noeth{\'e}rien, 
d{\'e}finisssons le $(n+1)$-foncteur:
\[
\begin{array}{rlcl}
I_{X/S}\; :& PIC^{n+1}(X) &\F &PIC(S)\\
& (L_1, \cdots , L_{n+1}) &\mapsto &
{\displaystyle \bigotimes_{k=0}^{n+1}
\left(  
\bigotimes_{i_1<\cdots <i_k} \gd (L_{i_1} \otimes \cdots \otimes L_{i_k})
\right)
^{(-1)^{n+1-k}}}
\end{array}\; ,
\]
appel{\'e} {\em foncteur fibr{\'e} d'intersection} pour $X/S$.
\begin{prop}
\label{propr-inter}
Le foncteur $I_{X/S}$ est un $(n+1)$-foncteur additif en chaque
variable, v{\'e}rifiant les propri{\'e}t{\'e}s suivantes:
\begin{enumerate}
\item La formation  de $I_{X/S}$ commute aux changements de bases.
\item $I_{X/S}$ est muni de donn{\'e}es de sym{\'e}trie, compatibles avec les
  donn{\'e}es d'additivit{\'e} en chaque variable.
\item Si $\pi : X\F S$ est fini et plat, le foncteur 
$I_{X/S}:PIC(X)\F PIC(S) $ est simplement le foncteur norme $N_{X/S}$
et les contraintes d'additivit{\'e} pour $I_{X/S}$ sont les isomorphismes
  usuels $N_{X/S} (L\otimes L\pr ) \isom N_{X/S} (L)
\otimes N_{X/S} (L\pr )$.
\item Soient $L_1,\cdots ,L_{n+1}$  des faisceaux inversibles sur
  $X$ et $\gs_{n+1}$ une section $\pi$-r{\'e}guli{\`e}re de $L_{n+1}$
  d{\'e}finissant un  diviseur de Cartier relatif effectif $D$. 
Il existe un isomorphisme canonique:
\[
\rho _D:
I_{X/S} (L_1, \cdots ,L_{n+1})
\isom 
I_{D/S} (L_1|_D, \cdots ,L_n|_D)
\]
qui est fonctoriel en les isomorphismes $L_i \isom L_i\pr$ pour $1\leq
i\leq n$ et
compatible avec les donn{\'e}es d'additivit{\'e}  et de sym{\'e}trie de 
$I_{X/S}$ 
en les $L_1,\cdots ,L_n$ et de $I_{D/S}$ en les 
 $L_1|_D,\cdots ,L_n|_D$.
\item Si, en plus des donn{\'e}es de 4, on dispose d'une section $\gs_n$
  de $L_n$ telle que $(\gs_n,\gs_{n+1})$ et $(\gs_{n+1},\gs_n)$ sont
  des suites $\pi$-r{\'e}guli{\`e}res, notons $E$ le diviseur d{\'e}fini par
  $\gs_n$. Le diagramme suivant est alors commutatif:
\[
\begin{CD}
I_{X/S} (L_1, \cdots ,L_{n+1}) @>{\rho_D}>> I_{D/S} (L_1|_D, \cdots ,L_n|_D)\\
@V{\rho_E}VV @VV{\rho_E}V\\
I_{E/S} (L_1|_E, \cdots ,L_{n-1}|_E,L_{n+1}|_E)
@>{\rho_D}>>
I_{D\cap E/S} (L_1|_{D\cap E}, \cdots ,L_{n-1}|_{D\cap E})
\end{CD}
\]
\end{enumerate}
\end{prop}
\begin{proof}[Preuve]
Notons que, par d{\'e}finition,
\[
I_{X/S}(L_1,\cdots ,L_n) =
\gd _{X/S}( K_{\ox}(L_1,\cdots ,L_n)) \; .
\]
C'est donc le $(n+1)$-foncteur multi-additif et sym{\'e}trique
associ{\'e} {\`a} la 
structure du cube sur le foncteur $\gd :PIC(X) \F PIC (S)$. 
La propri{\'e}t{\'e} 1 provient de la compatibilit{\'e} du foncteur $\gd$ aux
changements de bases. L'existence de donn{\'e}es d'additivit{\'e} en chaque
variable et de donn{\'e}es de sym{\'e}trie compatibles {\`a} l'additivit{\'e} 
provient de l'existence d'une structure de $(n+2)$-cube sur $\gd$.
 L'identification 3 de $I_{X/S}$ et $N_{X/S}$ dans
le cas d'un sch{\'e}ma relatif fini et plat est simplement la d{\'e}finition
de la structure du carr{\'e} dans le cas de dimension relative 0.\\
Il reste {\`a} construire l'isomorphisme  $\rho _D$: 
Les isomorphismes de restrictions 
d{\'e}crits en (\ref{th-prin})
induisent un isomorphisme canonique:
\[
\gd _{X/S} (K_{\ox} (\ox (D),L_1,\cdots ,L_n))
\isom 
\gd _{D/S} (K_{\ox (D)|_D}(L_1 (D)|_D, \cdots ,L_n (D)|_D))
\]
Par ailleurs, l'{\'e}galit{\'e}
\[
K_{\ox (D)|_D}(L_1(D)|_D,\cdots ,L_n(D)|_D)
=
K_{\oo _D}(L_1|_D,\cdots ,L_n|_D) \otimes_{\oo_D} \ox (D)|_D
\]
induit, d'apr{\`e}s (\ref{def-multifonct}), un isomorphisme
\[
\gd _{D/S} (K_{\ox (D)|_D}(L_1 (D)|_D, \cdots ,L_n (D)|_D))
\isom
I_{D/S} (L_1 |_D, \cdots ,L_n |_D).
\]
D'apr{\`e}s les propri{\'e}t{\'e}s des structures du cube pour $\gd _{D/S}$ et 
 $\gd _{X/S}$, 
ces deux isomorphismes sont fonctoriels et  additifs en chaque $L_i$
et compatibles aux donn{\'e}es de sym{\'e}tries. Notons $\rho _D$ leur
compos{\'e}; il v{\'e}rifie donc les propri{\'e}t{\'e}s demand{\'e}es.
\end{proof}
\begin{cor}
Soient $L_1,\cdots ,L_n$ des faisceaux inversibles sur $X$ et 
$(\gs_1,\cdots ,\gs_n)$ une suite $\pi$-r{\'e}guli{\`e}re de sections des
$L_i$. Notons $D_i$ le lieu des z{\'e}ros de $\gs_i$. Pour tout faisceau
inversible $L$ sur $X$, les
isomorphismes de restriction successifs de 
$\cap_{i=1}^p D_i$ {\`a} $\cap_{i=1}^{p+1} D_i$
induisent  un
isomorphisme canonique:
\[
I_{X/S} (L_1, \cdots ,L_n,L)
\isom 
N_{D_1\cap \cdots \cap D_n/S} (L|_{D_1\cap \cdots \cap D_n})\; ,
\]
fonctoriel en les isomorphismes $L\isom L\pr$ et additif en $L$.
\end{cor}
\begin{proof}[Preuve]
 Remarquons
que, par hypoth{\`e}se, pour tout entier $p$
tel que $1\leq p \leq n$, $D_1 \cap \cdots \cap D_p \F S$ est un morphisme
projectif et plat {\`a} fibres de dimension $n-p$ et 
$D_1 \cap \cdots \cap D_p\cap D_{p+1}$ est un diviseur de Cartier relatif
effectif de $D_1 \cap \cdots \cap D_p /S$.
On peut donc it{\'e}rer la construction pr{\'e}c{\'e}dente, et on obtient un
isomorphisme 
\[
I_{X/S} (\ox (D_1), \cdots , \ox (D_n),L)
\isom 
I_{D_1\cap \cdots \cap D_n/S} (L|_{D_1\cap \cdots \cap D_n})\; .
\]
On conclut alors en utilisant le fait que $D_1\cap \cdots \cap D_n \F
S$ est fini et plat, ce qui entra{\^\i}ne que
$I_{D_1\cap \cdots \cap D_n/S} (L|_{D_1\cap \cdots \cap D_n})
\simeq
N_{D_1\cap \cdots \cap D_n/S} (L|_{D_1\cap \cdots \cap D_n})$. La
fonctorialit{\'e} et l'additivit{\'e} de l'isomorphisme obtenu se
v{\'e}rifient {\`a} chaque {\'e}tape de l'it{\'e}ration.
\end{proof} 
\subsection{Sections du fibr{\'e} d'intersection}
\subsubsection{Cas de la dimension 0}
Soit $\pi :Y\F T$ est un morphisme fini et plat, et soit $L$ un
faisceau inversible sur $Y$, muni d'une section $\gs$ qui ne s'annulle
pas au dessus d'un ouvert $U=T\setminus Z$ de $T$. La section $\gs$
induit un isomorphisme $\oo |_{Y_U} \isom L|_{Y_U}$, qui induit donc
par fonctorialit{\'e} un isomorphisme 
$N_{Y_U/U} (\gs ):\oo _U \isom (N_{Y/T} L)|_U$. Celui ci se prolonge en
une section $s:\oo _T \isom N_{Y/T} L$, non nulle
en dehors de $Z$ et d{\'e}finie par le morphisme
\[
\det \pi \lst s: \det \pi \lst \oo _Y \F \det \pi \lst L \; .
\]
\subsubsection{Construction d'une section sur un ouvert}
Soit $\pi : X\F S$ un morphisme projectif et plat,
{\`a} fibres de dimension $n$ et $L_1,\cdots ,L_{n+1}$ des faisceaux
inversibles sur $X$.
Supposons que l'on dispose de sections $\gs_i$ de chaque $L_i$ telles
que $(\gs_1,\cdots ,\gs_n)$ est une suite $\pi$-r{\'e}guli{\`e}re. Notons
$D_i$ le  lieu des z{\'e}ros de $\gs_i$ et  $U$
l'ouvert de $S$ au dessus duquel  
  $D_1\cap \cdots \cap D_{n+1} = \emptyset$ , on peut appliquer ce
qui pr{\'e}c{\`e}de au morphisme fini et plat  $D_1\cap \cdots \cap D_n \F S$,
au faisceau inversible $L_{n+1}|_{D_1\cap \cdots \cap D_n}$
et {\`a} sa section $\gs_{n+1}|_{D_1\cap \cdots \cap D_n}$. On obtient
ainsi une section $\os \F I_{X/S}(L_1,\cdots ,L_{n+1})$, 
dont la restriction {\`a} $U$ est un isomorphisme.
\subsection{Construction du r{\'e}sultant.}
Soit $\pi : X\F S$ un morphisme projectif et plat,
{\`a} fibres de dimension $n$ et $L_1,\cdots ,L_{n+1}$ des faisceaux
inversibles sur $X$. Effectuons le changement de base:
\[
\xymatrix{X_P \ar[r]^-g\ar[d]_-{\pi} &X \ar[d]^-{\pi}\\
P= P_{L_1}\times _S \cdots \times _S P_{L_n} \ar[r]^-f & S
}
\]
Pour tout entier $i$, notons $p_i$ la projection de $P$ sur $P_i$ et
$\pi_i$ la compos{\'e}e $X\F P\F P_i$, et
consid{\'e}rons le faisceau inversible 
$L\pr_i=g\ust L_i \otimes \pi_i\ust \oo_{P_i}(1)$. Consid{\'e}rons alors le
faisceau inversible sur $P$:
\[
Res (L_1,\cdots ,L_{n+1}) = I_{X_P/P}(L_1\pr,\cdots ,L_{n+1}\pr)
\]
\begin{lemme}
\label{isom-resultant}
$Res (L_1,\cdots ,L_{n+1})$ est canoniquement isomorphe {\`a} 
$
f\ust I_{X/S}(L_1,\cdots ,L_{n+1}) \otimes 
{\displaystyle \bigotimes_{i=1}^{n+1}p_i\ust \oo_{P_i}(k_i)}
$, o{\`u} $k_i:S\F\Z$ est la fonction localement constante qui {\`a} $s\in S$
associe le nombre d'intersection sur $X_s$ des restrictions {\`a} $X_s$
des faisceaux $L_j$ pour $j\neq i$.
\end{lemme} 
\begin{proof}[Preuve]
Rappelons que, si $K$ est un $p$-cube dans $PIC(X)$, pour tout faisceau
inversible $L$ sur $S$, on a un isomorphisme canonique
$\gd(K\otimes\pi\ust L) \isom \gd(K)\otimes L^{\otimes\chi_{X/S}(K)}$. 
Soit $1\leq i \leq n+1$, notons $i_1,\cdots,i_n$  les $n$ entiers
compris entre 1 et $n+1$ et distincts de $i$, class{\'e}s par ordre
croissant. Si $M_{i_1},\cdots,M_{i_n}$ sont des faisceaux inversibles
sur $P$, on a alors:
\begin{multline*}
\chi_{X_P/P}(K_{\ox}(g\ust L_{i_1}\otimes \pi\ust M_{i_1} ,
\cdots ,g\ust L_{i_n}\otimes \pi\ust M_{i_n})=\\
\chi_{X_P/P}(K_{\ox}(g\ust L_{i_1},\cdots ,g\ust L_{i_n})=
\chi_{X/S}(K_{\ox}(L_{i_1},\cdots ,L_{i_n}))=k_i
\end{multline*}
 On peut alors {\'e}crire:
\begin{multline*}
Res (L_1,\cdots ,L_{n+1})
=
I_{X_P/P}(L\pr_1,\cdots ,L\pr_n,
g\ust L_{n+1}\otimes \pi_i\ust \oo_{P_{n+1}}(1))\\
\simeq
\gd(K)^{-1}\otimes \gd(K\otimes 
g\ust L_{n+1}\otimes \pi_i\ust \oo_{P_{n+1}}(1))
\simeq
\gd(K)^{-1}\otimes \gd(K\otimes 
g\ust L_{n+1}) \otimes( p_i\ust  \oo_{P_{n+1}}(1))^{k_{n+1}}\\
\simeq\gd(
\xymatrix{K\ar@{-}[r]_-{n+1}&K\otimes 
g\ust L_{n+1}})\otimes p_i\ust  \oo_{P_{n+1}}(k_{n+1}),
\end{multline*}
o{\`u} l'on a pos{\'e} $K=K_{\ox}(L_1\pr,\cdots,L_n\pr)$. 
En effectuant la m{\^e}me op{\'e}ration sur chacun des $L_i$, on obtient
finalement:
\begin{align*}
Res (L_1,\cdots ,L_{n+1})&
\simeq
I_{X_P/P}(g\ust L_1,\cdots,g\ust L_{n+1}) 
\otimes
\bigotimes_{i=1}^{n+1} p_i\ust  \oo_{P_i}(k_i)\\
&\simeq f\ust I_{X/S}(L_1,\cdots,L_{n+1}) 
\otimes
\bigotimes_{i=1}^{n+1} p_i\ust  \oo_{P_i}(k_i)
\end{align*}
\end{proof}
Pour tout $i$, consid{\'e}rons la section canonique $\gs_i$ des $L_i$ sur
$X_P$ et notons $D_i$ le lieu des z{\'e}ros de $\gs_i$.
On va construire une section canonique 
$Res(\gs_1,\cdots ,\gs_{n+1})$ de 
$Res(L_1\pr, \cdots , L_{n+1}\pr)$, qu'on appellera 
{\em r{\'e}sultant} des $\gs_i$ .\\
Pour toute permutation $\phi\in S_{n+1}$, notons $U_{\phi}$ l'ouvert
de $P$ au dessus duquel  $(\gs _{\phi (1)},\cdots ,\gs _{\phi (n)})$
est une suite $\pi$-r{\'e}guli{\`e}re.
\begin{lemme}
\label{genericite-suite}
Soient $U=\bigcup_{\phi\in S_{n+1}}U_{\phi}$ et 
$V=\bigcap_{\phi\in S_{n+1}}U_{\phi}$ et soient $Z$ et $Z\pr$ les
ferm{\'e}s compl{\'e}mentaires. Si $f:P\F S$ d{\'e}signe le morphisme
structural, on a:
\begin{enumerate}
\item 
Pour tout $z\in Z$, on a: $\text{Prof} (\oo_{f^{-1}(f(z)),z}) \geq 2$.
\item
Pour tout $z\in Z\pr$, on a: $\text{Prof}(\oo_{f^{-1}(f(z)),z})\geq
1$.
\end{enumerate}
\end{lemme}
\begin{proof}[Preuve]
Soit  $z\in Z$, soit $k<n$ la longueur maximale d'une suite
$\pi$-r{\'e}guli{\`e}re prise parmi les $\gs_i$ au dessus de $z$, soit
 $(\gs_{\phi (1)},\cdots ,\gs_{\phi (k)})$ une telle suite et soit
 $P\pr= P_{\phi (1)}\times _S\cdots\times _SP_{\phi (k)}$. D{\'e}composons  
$f:P\F S$ en 
\[
\xymatrix{
&P\pr\times _SP_{\phi (k+1)}\ar[rd]^-{q\pr} \\
P\ar[rr]^-p \ar[ru]^-{p\pr} \ar[rd]^-{p\sec}
&&P\pr \ar[r]^-r&S\\
&P\pr\times _SP_{\phi (k+2)}\ar[ru]^-{q\sec}
}
\]
Par le lemme (\ref{genericite2}), $p\pr (z)$ est contenu dans un
hyperplan de la fibre $(q\pr )^{-1} (p(z))$ et, de m{\^e}me,  
 $p\sec (z)$ est contenu dans un
hyperplan de la fibre $(q\sec )^{-1} (p(z))$. Donc $z$ est contenu
dans l'intersection de deux hyperplans de la fibre $p ^{-1}(p(z))$ et
donc dans l'intersection de deux hyperplans de $f ^{-1}(f(z))$, ce qui
montre le r{\'e}sultat.\\
La deuxi{\`e}me assertion se montre de mani{\`e}re analogue.
\end{proof}
Pour tout $\phi\in S_{n+1}$, on a un isomorphisme d{\'e}fini sur l'ouvert
$U_{\phi}$:
\begin{equation}
\label{isom-recur}
I_{X/S} (L_1,\cdots ,L_{n+1})
\isom
N_{D_{\phi(1)}\cap\cdots\cap D_{\phi(n)}/S} (L_{\phi(n+1)})
\end{equation}
et la section 
$N_{D_{\phi(1)}\cap\cdots\cap D_{\phi(n)}/S} (\gs_{\phi(n+1)})$ 
d{\'e}finit
donc une section $s_\phi$ de $I_{X/S}|_{U_\phi}$.
\begin{lemme}
Pour tous $\phi,\psi\in S_n$, on a  
 $s_\phi|_V=s_\psi|_V$.
\end{lemme}
\begin{proof}[Preuve]
 Si $\phi(n+1)=\psi(n+1)$, ceci
 provient de l'assertion 5 de la proposition (\ref{propr-inter}). 
Il suffit donc de montrer l'assertion 
dans le cas o{\`u} $\phi(n)=\psi(n+1)$ et $\phi(n+1)=\psi(n)$. 
Notons alors $Y=D_{\phi(1)}\cap\cdots\cap D_{\phi(n-1)}$. C'est un
 sch{\'e}ma relatif sur $S$ {\`a} fibres de dimension 1 et
les isomorphismes (\ref{isom-recur}) correspondant {\`a} $\phi$ et $\psi$
 se factorisent  en
\[
I_{X/S} (L_1,\cdots ,L_{n+1})
\isom
I_{Y/S}(L_{\phi(n)},L_{\phi(n+1)}) 
\isom
N_{Y\cap D_{\phi(n)}/S} (L_{\phi(n+1)})
\]
et 
\[
I_{X/S} (L_1,\cdots ,L_{n+1})
\isom
I_{Y/S}(L_{\phi(n)},L_{\phi(n+1)}) 
\isom
N_{Y\cap D_{\phi(n+1)}/S} (L_{\phi(n)})
\]
Le lemme (\ref{indep2}) entraine que le diagramme 
\[
\xymatrix{
&N_{Y\cap D_{\phi(n)}/S} (L_{\phi(n+1)})\ar[r]^-{\sim}&
I_{Y/S}(L_{\phi(n)},L_{\phi(n+1)})\ar[dd] \\ 
{\os} \ar[ur]^-{N_{Y\cap D_{\phi(n)}/S} (\gs_{\phi(n+1)})}
\ar[dr]_-{N_{Y\cap D_{\phi(n+1)}/S} (L_{\phi(n)})} \\
&N_{Y\cap D_{\phi(n+1)}/S} (L_{\phi(n)})\ar[r]^-{\sim}&
I_{Y/S}(L_{\phi(n)},L_{\phi(n+1)}) 
}
\]
est commutatif, ce qui montre l'assertion.
\end{proof}
En utilisant l'assertion 2 du lemme (\ref{genericite-suite}), on en
d{\'e}duit que $s_\phi$ et $s_\psi$ coincident sur $U_\phi\cap U_\psi$.
 Les $s_\phi$ d{\'e}finissent donc une section de $I_{X_P/P}$
sur l'ouvert $U$ de $P$. D'apr{\`e}s l'assertion 1 du lemme
(\ref{genericite-suite}), une telle section
se prolonge de mani{\`e}re unique en une section de $I_{X_P/P}$ sur
$P$.
\begin{Def}
On appelle r{\'e}sultant de $\gs_1,\cdots ,\gs_{n+1}$ et on note 
$Res(\gs_1,\cdots ,\gs_{n+1})$ l'uni\-que  section de
 $I_{X_P/P}$ qui coincide avec  $s_\phi$ sur l'ouvert $U_\phi$ pour
 tout $\phi \in S_{n+1}$.
\end{Def}
 Les r{\'e}sultats suivants justifient la terminologie "r{\'e}sultant":
\begin{lemme}
Soit  $R\subset P$ l'image par $\pi$ du lieu des z{\'e}ros de la section
$\gs =\bigoplus\gs_i$ de $\bigoplus L\pr_i$. $R$ est un ferm{\'e} de $P$ et 
pour tout point $p\in P$ tel que $\text{Prof}(\oo_{P_{f(p)},p})\leq 1$,
$p$ est {\'e}l{\'e}ment de $R$ si et seulement si $Res(\gs_1,\cdots
,\gs_{n+1})(p)=0$. 
\end{lemme}
\begin{proof}[Preuve]
 D'apr{\`e}s la partie 2 du lemme (\ref{genericite-suite}), il suffit de
 montrer que, pour tout point $p$ de l'ouvert 
$U=\bigcup_{\phi\in S_{n+1}}U_{\phi}$, on a
 $Res(\gs_1,\cdots ,\gs_{n+1})(p)=0$ si et
 seulement si $p\in R$.  Soit donc  $p\in U$ et $\phi\in S_{n+1}$ tel
que $p\in U_\phi$. Sur $U_\phi$, $Res(\gs_1,\cdots ,\gs_{n+1})$
 coincide avec 
$N_{D_{\phi(1)}\cap\cdots\cap D_{\phi(n)}} (\gs_{\phi(n+1)})$ et cette
 section s'annulle si et seulement si $\gs_{\phi(n+1)}$ s'annulle en
 point de la fibre $(D_{\phi(1)}\cap\cdots\cap D_{\phi(n)})_p$, ce qui
 montre le r{\'e}sultat.
\end{proof}

\begin{lemme}
\label{sym-resultant}
Pour toute permutation $\phi\in S_{n+1}$, l'isomorphisme
 de sym{\'e}trie 
\[
Res(L_{\phi(1)},\cdots,L_{\phi(n+1)})
\isom
Res(L_1,\cdots,L_{n+1})
,
\]
 provenant des isomorphismes de sym{\'e}trie
du fibr{\'e} d'intersection, {\'e}change leurs sections respectives 
$Res(\gs_{\phi(1)},\cdots,\gs_{\phi(n+1)})$
et
$Res(\gs_1,\cdots,\gs_{n+1})$.
\end{lemme}
\begin{proof}[Preuve]
Il suffit de montrer cette {\'e}galit{\'e} sur l'ouvert 
$V=\bigcap_{\phi\in S_{n+1}}U_{\phi}$ et de plus il suffit
de consid{\'e}rer le cas o{\`u} $\phi$ est une transposition de deux indices
$i$ et $j$. 
Si $n=1$, il suffit alors d'appliquer le lemme (\ref{indep2}). Si
$n>1$, choisissons une permutation $\psi\in S_n$ telle que 
$\psi(n+1)\notin \{ i,j\} $. Comme, sur $V$, la section 
$Res(\gs_1,\cdots,\gs_{n+1})$ se factorise en:
\[
\os \F N_{D_{\psi(1)}\cap\cdots\cap D_{\psi(n)}/S} (L_{\psi(n+1)})
\isom  I_{X/S}(L_1,\cdots,L_{n+1})\; ,
\]
le r{\'e}sultat en r{\'e}sulte, puisque $D_{\psi(1)}\cap\cdots\cap
D_{\psi(n)}$ est invariant par la transposition $\phi$.
\end{proof}

\addtocounter{subsubsection}{6}
\subsubsection{Multiplicativit{\'e} du r{\'e}sultant}
\addtocounter{theo}{1}
Soient $L_1\pr,L_1\sec,L_2,\cdots,L_{n+1}$ des faisceaux inversibles
suf\-fisamment positifs sur $X$ et supposons de plus que 
$L_1\pr\otimes L_1\sec\gg 0$. On notera alors 
$P\pr=P_{L_1\pr}\times_S\cdots\times_SP_{L_{n+1}}$, 
$P\sec =P_{L_1\sec}\times_S\cdots\times_SP_{L_{n+1}}$, $P=P\pr \times_S
P\sec$ et  $Q=P_{L_1\pr\otimes L_1\sec}\times_S\cdots\times_SP_{L_{n+1}}$. Le 
 morphisme canonique de multiplication des sections 
$m:P_{L_1\pr}\times_SP_{L_1\sec}\F P_{L_1\pr\otimes L_1\sec}$ induit un
morphisme not{\'e} toujours $m: P\F Q$.
Consid{\'e}rons alors le diagramme suivant
\begin{equation}
\label{diagramme-resultant}
\xymatrix{
&X_P\ar[d]_-{\pi}\ar[rr]_-n \ar@/_/[ddl]_-{\pi\pr}
 \ar@/^/[ddr]^-{\pi\sec} \ar@/^1pc/[rrr]^g
&&X_Q\ar[d]_-{\pi}\ar[r]_-f\ar@/^1pc/[ddd]^{\pi_1}
&X\ar[d]_-{\pi} \\
&P \ar[rr]_-m \ar[dl]^-{p\pr}\ar[dr]_-{p\sec}
&&Q\ar[r] \ar[dd]_-{p_1}&S\\
P\pr \ar[d]_-{p_1\pr}&&P\sec \ar[d]_-{p_1\sec}\\
P_{L_1\pr} &&P_{L_1\sec}& P_{L_1\pr\otimes L_1\sec}
}
\end{equation}
On notera enfin  $p_i$ les diff{\'e}rentes
projections de $P$, $P\pr$ ou $Q$ sur $P_{L_i}$ et $\pi_i$ leur
compos{\'e}e avec $\pi$. 
\begin{lemme} 
\label{mult-resultant}
Il existe  un isomorphisme
canonique de faisceaux inversibles sur  $P\pr\times_SP\sec$
\[
m\ust Res(L_1\pr\otimes L_1\sec,L_2,\cdots,L_{n+1})
\isom
(p_1\pr)\ust  Res(L_1\pr,L_2,\cdots,L_{n+1})
\otimes
(p_1\sec)\ust  Res(L_1\sec,L_2,\cdots,L_{n+1})
\]
qui induit une {\'e}galit{\'e} sur les sections: 
\[
m\ust Res(\tau,\gs_2,\cdots,\gs_{n+1})
=
(p_1\pr)\ust  Res(\gs_1\sec,\gs_2,\cdots,\gs_{n+1})
\otimes
(p_1\sec)\ust  Res(\gs_1\pr,\gs_2,\cdots,\gs_{n+1})
\]
o{\`u} $\tau$ est la section canonique de 
$f\ust (L_1\otimes L_1\pr)\otimes 
\pi\ust_1 \oo_{P_{L_1\pr \otimes L_1\sec}}(1)$.
\end{lemme}
\begin{proof}[Preuve]
Il existe un isomorphisme canonique de faisceau inversibles sur 
 $P\pr\times_SP\sec$:
\[
m\ust p_1\ust\oo_{P_{L_1\pr\otimes L_1\sec}}(1) 
\simeq 
(p_1\pr)\ust\oo_{P_{L_1\pr}}(1) 
\otimes
(p_1\sec)\ust\oo_{P_{L_1\sec}}
\]
En prenant l'image inverse par $\pi$ et en  tensorisant avec 
$g\ust(L_1\pr\otimes L_1\sec)$, on obtient sur
$X_{P\pr\times_SP\sec}$ un isomorphisme:
\[
n\ust \pi\ust\oo_{P_{L\otimes M}}(1) \otimes g\ust(L_1\pr\otimes L_1\sec)
\simeq 
\left(
(\pi_1\pr)\ust \oo_{P_{L_1\pr}}(1) \otimes g\ust L_1\pr
\right)
\otimes
\left(
(\pi_1\sec)\ust \oo_{P_{L_1\sec}}(1) \otimes g\ust L_1\sec
\right)
\]
qui identifie $n\ust\tau$ avec $\gs_1\pr \otimes \gs_1\sec$.
La prori{\'e}t{\'e} de  multiplicativit{\'e} du fibr{\'e}
d'intersection $I_{X_Q/Q}$ 
induit donc un isomorphisme:
\begin{multline*}
m\ust I_{X_Q/Q}
(f\ust (L_1\pr\otimes  L_1\sec) \otimes \pi_1\ust
\oo_{P_{L_1\pr\otimes L_1\sec}}(1),\cdots, 
f\ust L_{n+1}\otimes \pi_n\ust \oo_{P_{L_{n+1}}}(1))
\isom \\
 I_{X_Q/Q}(g\ust L_1\pr\otimes g\ust L_1\sec \otimes (\pi_1\pr)\ust
\oo_{P_{L_1\pr}}(1)
\otimes (\pi_1\sec)\ust
\oo_{P_{L_1\sec}}(1)
,\cdots, 
g\ust L_{n+1}\otimes \pi_n\ust \oo_{P_{L_{n+1}}}(1))
\\
\isom 
(p_1\pr)\ust  Res(L_1\pr,\cdots,L_{n+1})
\otimes
(p_1\sec)\ust  Res(L_1\sec,\cdots,L_{n+1})
\end{multline*}
qui identifie les sections correspondantes.
\end{proof}
\subsection{Lien avec la d{\'e}finition classique du r{\'e}sultant}
Consid{\'e}rons le sch{\'e}ma $X=\mathbb{P}^n_k$ sur $S=Spec(k)$. Dans tout ce
qui suit, on {\'e}crira explicitement $X$ sous la forme
$Proj(k[X_1,\cdots,X_{n+1}]$, et on consid{\`e}rera les  $n+1$ hyperplans
$H_i$ de $X$ d'{\'e}quations $X_i=0$. Fixons des
entiers strictement positifs $d_1,\cdots,d_{n+1}$ et consid{\'e}rons les 
faisceaux $L_i=\ox(d_i)$, qu'on identifiera {\`a} $\ox(d_iH_i)$, en notant
$\tau_i$ la section canonique correspondante. Les $L_i$
 sont suffisamment positifs et on regarde les espaces
projectifs $P_i$ correspondants sont de dimension
$\binom{d_i+n}{n}-1$. Notons enfin $V_i=H^0(X,L_i)$.\\ 
L'isomorphisme canonique (\ref{isom-resultant}) s'{\'e}crit ici: 
\[
Res (L_1,\cdots,L_{n+1})
\isom
\bigotimes_{i=1}^{n+1} p_i\ust \oo_{P_i}(k_i)
\otimes
I_{X/S}(L_1,\cdots,L_{n+1})
\]
avec $k_i=\prod_{k\neq i}d_k$. Comme
 $\bigcap_{i=1}^{n+1}H_i=\emptyset$, $I_{X/S}(\tau_1,\cdots,\tau_{n+1})$
 d{\'e}finit une trivialisation de $ I_{X/S}(L_1,\cdots,L_{n+1})$. On en
 d{\'e}duit donc un 
isomorphisme, associ{\'e} canoniquement au choix des coordonn{\'e}es $X_i$:
\[
Res (L_1,\cdots,L_{n+1})
\isom
\bigotimes_{i=1}^{n+1} p_i\ust \oo_{P_i}(k_i)
\]
La section canonique 
$Res (\gs_1,\cdots,\gs_{n+1})$ de $Res (L_1,\cdots,L_{n+1})$ 
d{\'e}finit donc une section de 
$\bigotimes_{i=1}^{n+1} p_i\ust \oo_{P_i}(k_i)$
 sur le produit 
$P=P_1\times\cdots\times P_{n+1}$, c'est {\`a} dire un polyn{\^o}me
quasi-homog{\`e}ne sur $V=\bigoplus_{i=1}^{n+1}V_i$, de degr{\'e}
$k_i$ relativement {\`a} $V_i$. C'est le r{\'e}sultant  {\'e}tudi{\'e} par
\textsc{Jouanolou} dans \cite{jouanolou1}. On le notera alors 
$\underline{Res}(\gs_1,\cdots,\gs_{n+1})$. Notons alors que les lemmes
(\ref{sym-resultant}) et (\ref{mult-resultant}) entrainent que ce
polyn{\^o}me est sym{\'e}trique etmultiplicatif en chaque groupe de variables.
\subsubsection{La formule de Poisson}
Dans cette section, on interpr{\`e}te la formule de Poisson pour les
polyn{\^o}mes r{\'e}sultants (\cite{jouanolou1}, Prop 2.7,p124), qui permet de
calculer le 
r{\'e}sultant $\underline{Res}(\gs_1,\cdots,\gs_{n+1})$ par r{\'e}currence sur
la dimension, en l'exprimant en fonction 
d'un r{\'e}sultant $\underline{Res}(\gs_1|_H,\cdots,\gs_n|_H)$
associ{\'e} {\`a} un 
hyperplan $H$, et de la norme d'une fonction 
d{\'e}finie sur l'ouvert affine, compl{\'e}mentaire de $H$.\\
On note ici $H$ l'hyperplan $H_{n+1}$
de $X$ et $Q=P_1\times\cdots\times P_n$ de telle sorte que $P=Q\times
P_{n+1}$. Soit $Z\subset X_P$ le lieu des z{\'e}ros communs de
$\gs_1,\cdots ,\gs_n$, consid{\`e}rons alors les ouverts 
$U_0=\{ q\in Q| Z\cap H_{n+1}\cap X_q=\emptyset\}$ et 
$U=U_0\times P_{n+1}\subset P=Q\times P_{n+1}$. $H$ est un espace
projectif de dimension $n-1$, muni de coordonn{\'e}es $X_1,\cdots,X_n$, et
on peut donc d{\'e}finir un r{\'e}sultant 
$Res(L_1|_H,\cdots,L_n|_H)$ sur $P_1\pr\times\cdots P_n\pr$, o{\`u}
$P_i\pr=\mathbb{P} (H^0(H,L_i)\ve)$.
R{\'e}sumons d'abord dans un diagramme la situation:
\[
\xymatrix{&&Res(\gs_1,\cdots,\gs_{n+1}) \ar@{.}[dr]
&X_U\subset X_P \ar[d]
 \\
Res(\gs_1|_H,\cdots,\gs_n|_H)\ar@{.}[dr] &
H_{Q\pr}\ar[d] &X_{U_0}\subset X_Q \ar[d] 
& P \ar[dl]_-{\pi_Q}\ar[dr]_-{\pi_{n+1}}
& X_{P_{n+1}}\ar[d]\ar[r]_-p &X \\
&Q\pr \ar[d]_-{\pi\pr_i} 
&U_0\subset Q\ar[l]_-{\pi_{Q\pr}}\ar[d]_-{\pi_i} &&P_{n+1} \\
&P_i\pr &P_i
}
\]
Sur l'ouvert $X\setminus H$, on dispose d'un isomorphisme 
$L_{n+1} \isom \oo_{X\setminus H}$. La restriction $\Tilde{\gs}_{n+1}$
de $\gs_{n+1}$ {\`a}  
$P\times (X\setminus H)$ est donc une section de
$p_n\ust(\oo(1))$. Comme $Z_U$ est contenu dans
$(X\setminus H)\times U$, la restriction de $\Tilde{\gs}_{n+1}$ {\`a}
$Z_U$ est bien d{\'e}finie et sa norme $N_{Z_U/U}(\Tilde{\gs}_{n+1})$ est
une section de  
$\pi_{n+1}\ust \oo(k_{n+1})$. Notons que ces consid{\'e}rations montrent
l'existence d'un isomorphisme canonique 
\begin{equation}
\label{poisson}
Res(L_1,\cdots,L_{n+1})|_U \isom 
\pi_{n+1}\ust \oo_{P_{n+1}}(k_{n+1})
\end{equation}
 qui identifie les sections $Res(\gs_1,\cdots,\gs_{n+1}|_U$ et 
 $N_{Z_U/U}(\Tilde{\gs}_{n+1})$.
\begin{prop}[Formule de Poisson]
Sur l'ouvert $U$, on a l'{\'e}galit{\'e} entre polyn{\^o}mes de 
 $\Gamma(U_0,\oo_{U_0})[V_{n+1}]$ :
\[
\underline{Res}(\gs_1,\cdots,\gs_{n+1})|_U
=
\pi_{Q\pr}\ust \underline{Res}(\gs_1|_H,\cdots,\gs_n|_H)|_{U_0}
.
N_{Z_U/U}(\Tilde{\gs}_{n+1})
\]
\end{prop}
\begin{proof}[Preuve]
Utilisons la d{\'e}finition du faisceau r{\'e}sultant et la
multiplicativit{\'e}  
du fibr{\'e} d'inter\-sec\-tion pour {\'e}crire le diagramme d'isomorphismes:
\[
\begin{array}{ccccc}
Res(L_1,\cdots,L_{n+1}) &
\isom &
I_{X/S}(L\pr_1,\cdots,L_n\pr,p\ust \ox(d_{n+1}H))&
\otimes &
I_{X/S}(L\pr_1,\cdots,L_n\pr,\pi_{n+1}\ust\oo_{P_{n+1}}(1))\\
\| &&\downarrow\wr && \downarrow\wr\\
Res(L_1,\cdots,L_{n+1})&&
\left( I_{H/S}(L\pr_1|_H,\cdots,L\pr_n|_H)\right)^{d_{n+1}} &&
N_{Z/P}(\pi_{n+1}\ust\oo_{P_{n+1}}(1))\\
\downarrow\wr&&\downarrow\wr && \downarrow\wr\\
\begin{array}{c}
\bigotimes_{i=1}^{n+1} \pi_i\ust\oo_{P_i}(k_i)\\
\otimes\\
I_{X/S}(L_1,\cdots,L_{n+1})
\end{array} &&
\begin{array}{c}
\bigotimes_{i=1}^n \pi_i\ust\oo_{P_i}(k_i)\\
\otimes \\
\left( I_{H/S}(L_1|_H,\cdots,L_n|_H)\right)^{d_{n+1}}
\end{array} && 
\begin{array}{ccc}
\pi_{n+1}\ust \oo_{P_{n+1}}(k_{n+1})\\
{}\\
{}
\end{array}
\end{array}
\]
La deuxi{\`e}me ligne de ce diagramme donne un isomorphisme
\begin{equation}
\label{poisson2}
Res(L_1,\cdots,L_{n+1})\isom 
\left( I_{H/S}(L\pr_1|_H,\cdots,L\pr_n|_H)\right)^{d_{n+1}} \otimes
N_{Z/P}(\pi_{n+1}\ust\oo_{P_{n+1}}(1))
\end{equation}
En combinant  la restriction de (\ref{poisson2}) {\`a} $U$ et
l'isomorphisme 
$Res(\gs_1|_H,\cdots,\gs_n|_H)|_U:
\oo_U\isom I_{H/S}(L\pr_1|_H,\cdots,L\pr_n|_H)|_U$, on obtient
l'isomorphisme (\ref{poisson}).
La restriction de (\ref{poisson2}) {\`a} $U$ identifie donc 
$Res(\gs_1,\cdots,\gs_{n+1})|_U$ {\`a}
=
$\pi_{Q\pr}\ust Res(\gs_1|_H,\cdots,\gs_n|_H)|_{U_0}
\otimes 
N_{Z_U/U}(\Tilde{\gs}_{n+1})$. L'{\'e}galit{\'e} polyn{\^o}miale
recherch{\'e}e s'en 
d{\'e}duit en utilisant la d{\'e}finition des polyn{\^o}mes
r{\'e}sultants et les 
isomorphismes de la
partie inf{\'e}rieure du diagramme.
\end{proof}

\end{document}